\documentclass[sn-mathphys,Numbered]{sn-jnl}
\usepackage{amsmath,amsfonts}
\newtheorem{remark}{Remark}
\newtheorem{lemma}{Lemma}
\usepackage{hyperref}
\usepackage{manyfoot}
\newcommand{\RR}{\mathbb{R}}
\newcommand{\NN}{\mathbb{N}}
\newcommand{\CC}{\mathbb{C}}
\newcommand{\rme}{\mathrm{e}}
\newcommand{\rmi}{\mathrm{i}}
\newcommand{\bb}{\boldsymbol}
\usepackage{pgfplots}
\pgfplotsset{compat=newest}
\usetikzlibrary{plotmarks}
\usetikzlibrary{arrows.meta}
\usepgfplotslibrary{patchplots}
\usepackage{grffile}

\begin{document}

\title{A $\mu$-mode approach for exponential
  integrators: actions of
  $\varphi$-functions of Kronecker sums}
\author*[1]{\fnm{Marco} \sur{Caliari}}\email{marco.caliari@univr.it}
\equalcont{These authors contributed equally to this work.}
\author[1]{\fnm{Fabio} \sur{Cassini}}\email{fabio.cassini@univr.it}
\equalcont{These authors contributed equally to this work.}
\author[2]{\fnm{Franco} \sur{Zivcovich}}\email{franco.zivcovich@gmail.com}
\equalcont{These authors contributed equally to this work.}

\affil[1]{\orgdiv{Department of Computer Science},
  \orgname{University of Verona}, \orgaddress{\street{Strada~Le~Grazie,~15},
    \city{Verona}, \postcode{37134}, \country{Italy}}}
\affil[2]{\orgname{Neurodec}, \orgaddress{\city{Sophia Antipolis}, \country{France}}}

\abstract{
  We present a method for computing actions of the exponential-like
  $\varphi$-functions
  for a Kronecker sum $K$ of $d$ arbitrary matrices $A_\mu$. It is based on
  the approximation of the integral
  representation of the $\varphi$-functions by Gaussian quadrature formulas
  combined with a scaling and squaring technique.
The resulting algorithm, which we call \textsc{phiks}, evaluates the required
actions
by means of $\mu$-mode products involving
exponentials of the \emph{small sized} matrices~$A_\mu$, without forming the
\emph{large sized} matrix~$K$ itself.
\textsc{phiks}, which profits from the highly efficient level 3 BLAS,
is designed to compute different $\varphi$-functions applied on the same
  vector or a linear combination of actions of \mbox{$\varphi$-func\-tions}
  applied on different vectors.
  In addition, thanks to
  the underlying scaling and squaring techniques,
  the desired quantities are available
  simultaneously at suitable time scales.
  All these features allow the effective usage of \textsc{phiks} in the
  exponential integration context.
  In fact, our newly designed  method has been tested
  in popular exponential
  Runge--Kutta integrators of stiff order from
  one to four, in comparison with state-of-the-art algorithms for computing
  actions of $\varphi$-functions. The numerical experiments with
  discretized semilinear evolutionary 2D or 3D
  advection--diffusion--reaction, Allen--Cahn, and Brusselator equations
  show the superiority of the proposed $\mu$-mode approach.}



\maketitle

\section{Introduction}
{\color{black}
We consider a system of Ordinary Differential Equations (ODEs) of the form
\begin{subequations}\label{eq:ODE}
\begin{equation}\label{eq:ODEinit}
  \left\{
  \begin{aligned}
    \bb u'(t)&=K\bb u(t)+\bb g(t,\bb u(t)),&t\in[0,T],\\
    \bb u(0)&=\bb u_0.
    \end{aligned}
  \right.
\end{equation}
Here, $\bb u\colon [0,T]\to\CC^N$ is the unknown,
being $N$ the total number of degrees of freedom,
$\bb g\colon [0,T]\times \CC^{N}\to\CC^N$ is a nonlinear function, and
$K\in\CC^{N\times N}$ is a
large sized matrix
which can be written as a Kronecker sum, that is
\begin{equation}\label{eq:kronsum}
  K=A_d\oplus A_{d-1}\oplus\cdots\oplus A_1=\sum_{\mu=1}^dA_{\otimes\mu},\quad
  A_{\otimes\mu}=I_d\otimes
  \cdots\otimes I_{\mu+1}\otimes A_\mu\otimes
  I_{\mu-1}\otimes \cdots \otimes I_1.
\end{equation}
\end{subequations}
Here and throughout the paper $d\in\NN$, $\mu=1,\ldots,d$,
$I_\mu$ is the identity matrix of size $n_\mu\times n_\mu$, the symbol $\otimes$
denotes the Kronecker product, and, unless
otherwise stated, $A_\mu\in\CC^{n_\mu\times n_\mu}$ is an arbitrary matrix.

These systems may arise
when applying the method of lines
to some evolutionary Partial Differential Equations (PDEs),
from different fields of science and engineering, defined
in a spatial domain $\Omega\subseteq \RR^d$
which is the Cartesian product of $d$ intervals.
Typical instances
are semilinear advection--diffusion, nonlinear Schr\"odinger,
or complex Ginzburg--Landau equations, possibly fractional in space.
In dimension $d=2$, for example,
the Laplace operator $\Delta=\partial_{x_1x_1}+\partial_{x_2x_2}$ on
a rectangular domain $\Omega$ with homogeneous Dirichlet boundary conditions
can be discretized  as
$K=A_2\oplus A_1=I_2\otimes A_1+A_2\otimes I_1\in\RR^{N\times N}$,
where $A_\mu\in\RR^{n_\mu\times n_\mu}$ is the discretization matrix
of $\partial_{x_\mu x_\mu}$ with standard second-order
finite differences and $N=n_1n_2$.
Notice that several other tensor-product approximation techniques
lead to Kronecker sums. We mention higher-order
(non)uniform finite differences and
lumped mass finite elements (that yield sparse matrices $A_\mu$),
or spectral differentiations (that yield dense matrices $A_\mu$).
Other type of boundary conditions can be considered,
as long as it is possible to write $K$ as a Kronecker sum. In particular,
inhomogeneous boundary conditions of Dirichlet or Neumann
type can be encapsulated into
the nonlinear term $\bb g$.}

When 
system~\eqref{eq:ODE}
is \emph{stiff}, a prominent way to numerically integrate it in time
is by using explicit exponential integrators~\cite{HO10}. These schemes
require the action of the exponential and/or of the
so-called $\varphi$-functions which, for a general matrix $X\in\CC^{N\times N}$,
are defined as
\begin{equation}\label{eq:phifunctions}
  \varphi_\ell(X)=\int_{0}^1f_{\ell}(\theta,X)d\theta,\quad
  f_\ell(\theta,X)=\frac{\theta^{\ell-1}}{(\ell-1)!}\exp((1-\theta)X),\quad
  \ell\ge 1.
\end{equation}
The direct approximation of the matrix $\varphi$-functions is
feasible only when the size of $X$
is not too large.
In this case, the most commonly employed  techniques are based on
Pad\'e
approximations~\cite{BSW07}, although
  other rational methods based on the numerical inversion
  by a quadrature formula of the Laplace
  transform \cite{LF10} or polynomial methods based on the truncated
  Taylor series \cite{LYL22} can be considered.
On the other hand, in this manuscript we are interested in
matrices of large size. In this case,
it is possible to rely on methods which directly compute the action
of the matrix $\varphi$-functions
on a vector, or even their linear combination at once.
Among them,
Krylov methods~\cite{GRT18,LPR19,NW12} and
other polynomial methods~\cite{AMH11,CCZ20,CCZ23,CKOR16,CKZ18,LYL22}
are typically employed.

{\color{black}When we consider matrices $K$ which are Kronecker
  sums~\eqref{eq:kronsum}, it is possible
to express the action of $\exp(K)$ on a vector $\bb v$ by
using the Kronecker product of the exponentials of the matrices
$A_\mu$. In fact, considering again
the two-dimensional case $K=A_2\oplus A_1$, we have
\begin{equation*}
  \begin{split}
\exp(K)&=\exp(I_2\otimes A_1+A_2\otimes I_1)=
\exp(I_2\otimes A_1)\exp(A_2\otimes I_1)\\
&=(I_2\otimes \exp(A_1))(\exp(A_2)\otimes I_1)=
\exp(A_2)\otimes\exp(A_1),
\end{split}
\end{equation*}
where we used the commutativity of
$I_2\otimes A_1$ and $A_2\otimes I_1$, the
semigroup property of the exponential function, and
the mixed-product property of the Kronecker product.
Therefore, we obtain
$\exp(K)\bb v=(\exp(A_2)\otimes\exp(A_1))\bb v$.}
The same result
can be accomplished 
by the more computationally attractive formula
$\exp(A_1)\bb V\exp(A_2)^{\sf T}$,
without computing
the Kronecker product $\exp(A_2)\otimes\exp(A_1)$
(see References~\cite{N69,BS17}).
Here, $\bb V$ is the matrix of size $n_1\times n_2$ whose
$j$-th column is
made by the elements of $\bb v$ from $(j-1)n_1+1$ to $jn_1$,
for $j=1,\ldots,n_2$. This approach
can be generalized to the
computation of the action of the matrix exponential when $K$ is the Kronecker
sum of $d$ \emph{arbitrary} matrices $A_\mu$
through the so-called
$\mu$-mode product {\color{black} (see the next section for more details)}.
With this technique, it is possible to efficiently implement exponential
schemes which require the action of the matrix exponential
only, such as some splitting schemes, Lawson methods and
Magnus integrators.
Unfortunately, such an elegant approach does not directly apply to the
computation of the action of $\varphi$-functions, since
they do not enjoy the aforementioned semigroup property of the exponential
function.

In this paper, we aim at computing
  actions of $\varphi$-functions for
a matrix $K$ which is the Kronecker sum of $d$
matrices $A_\mu$ using
a $\mu$-mode approach, without assembling the matrix $K$ itself.
Moreover, since we are interested in the application to exponential
integrators which may require more than a single $\varphi$-function
evaluation, as in
the case of
exponential Runge--Kutta schemes of high stiff order~\cite{HO05,L21},
we will derive an algorithm for the
computation of actions of
different $\varphi$-functions on the \emph{same}
vector
\begin{equation}\label{eq:diffphi}
  \{\exp(\tau K)\bb v,\varphi_1(\tau K)\bb v,\varphi_2(\tau K)\bb v,\ldots,
  \varphi_p(\tau K)\bb v\}
\end{equation}
at once, as well as for
\emph{linear combinations} of actions of $\varphi$-functions
\begin{equation}\label{eq:lincomb}
  \exp(\tau K)\bb v_0+\varphi_1(\tau K)\bb v_1+
  \varphi_2(\tau K)\bb v_2+
  \cdots+
  \varphi_p(\tau K)\bb v_p,
\end{equation}
where $\tau$ is the time step size of the integrator.
  For an efficient computation of the quantities in formula~\eqref{eq:diffphi}
  we will use the scaling and modified squaring method proposed in
  Reference~\cite{SW09}, while for the linear combination in
  formula~\eqref{eq:lincomb} we will derive a new scaling and squaring
  procedure.
  As a byproduct of these techniques, our algorithm can also output
  the desired quantities at different time scales, a feature of great
  importance   in the implementation of certain exponential integrators.

The remaining part of the paper is structured as follows.
In Section~\ref{sec:motivating}
we {\color{black} briefly describe the founding basis for the new method,
  \mbox{i.e.},
  the $\mu$-mode product, the Tucker operator, and its
  coupling with quadrature formulas for the evaluation of
  $\varphi$-functions.}
Section~\ref{sec:phiapprox}, the main one,
is devoted to the description of the new algorithm,
which we call \textsc{phiks} (PHI-functions of Kronecker Sums),
for the approximation of actions of $\varphi$-functions
on the same vector
and for linear combinations of actions of $\varphi$-functions.
An important subsection describes a suitable choice of the
  scaling parameter and of the quadrature formula. In particular, for the
  latter, we propose an effective  closed Gaussian formula.
Then, in Section~\ref{sec:numerical}, we validate our implementation
of \textsc{phiks} by running several
  examples  in dimensions $d=3$ and $d=6$ which involve different
$\varphi$-functions and their linear combinations. Moreover,
  we apply the proposed technique to the numerical solution of
  physically relevant 2D and 3D
  stiff advection--diffusion--reaction equations,
  with up to $N=2\cdot 100^3$ degrees of freedom,
and five different
exponential Runge--Kutta integrators of order up to four.
Finally, we draw some conclusions in Section~\ref{sec:conclusions}.
\section{A \texorpdfstring{$\mu$-mode}{mu-mode} approach for
  evolutionary equations in {K}ronecker form}\label{sec:motivating}
{\color{black}
  The founding basis of the technique 
  that we propose in this manuscript is based on
  the $\mu$-mode approach for the action of the matrix exponential.
  Due to its importance, we briefly recall here the main concepts,
  and invite a reader not familiar with the following formalism
to check References~\cite{CCEOZ22,CCZ23kp,KB09} for a thorough explanation
with full details. Let us denotes by $\bb V$ an order-$d$ tensor
of size $n_1\times\cdots \times n_d$ with elements $v_{i_1\ldots i_d}$,
and by $L_\mu$ a matrix of size $n_\mu\times n_\mu$ of elements
$\ell_{i j}^\mu$. Then, the
\emph{$\mu$-mode}
product of the tensor $\bb V$ with the matrix $L_\mu$, denoted
as $\bb V\times_\mu L_\mu$, is the tensor $\bb W$ of
size $n_1\times\cdots \times n_d$ defined elementwise as
\begin{equation*}
  w_{i_1\ldots i_d} = \sum_{j_\mu=1}^{n_\mu}
                      v_{i_1\ldots i_{\mu-1}j_\mu i_{\mu+1}\ldots i_d}
                      \ell_{i_\mu j_\mu}^\mu.
\end{equation*}
This corresponds to multiply the matrix
$L_\mu$ onto the \emph{$\mu$-fibers} of the tensor $\bb V$
(i.e., vectors along direction $\mu$ which are
generalizations to tensors of columns and rows of a matrix).
The concatenation of $\mu$-mode products with the matrices
$L_1,\ldots,L_d$, that is the tensor $\bb W$ with elements
\begin{equation*}
  w_{i_1\ldots i_d} = \sum_{j_d=1}^{n_d}\cdots
                      \sum_{j_1=1}^{n_1}v_{j_1\ldots j_d}
                      \prod_{\mu=1}^d \ell_{i_\mu j_\mu}^\mu,
\end{equation*}
is denoted by $\bb V\times_1 L_1\times_2\cdots\times_d L_d$ and
referred to as \emph{Tucker operator}. In terms of computational cost,  a
single $\mu$-mode product requires~$\mathcal{O}(Nn_\mu)$
floating point operations, being $N=n_1\cdots n_d$,
and it can be implemented by a single matrix-matrix product.
Consequently, the Tucker operator
has an overall computational cost of $\mathcal{O}(n^{d+1})$
for the case $n_1 = \ldots = n_d=n$.
It can be realized with $d$ calls of level 3 BLAS (Basic Linear
Algebra Subprograms \cite{DDCHD90}), whose
highly optimized implementations are available for any kind of modern
computer hardware (see, for instance, References~\cite{mkl,XQY12,cublas}).

The relation between the Kronecker product and the Tucker operator is given
by the following Lemma
(the proof can be found, for instance, in
Reference~\cite[Lemma 2.1]{CCZ23kp}).
\begin{lemma}\label{lem:krontommp}
Let $L_\mu \in \CC^{n_\mu\times n_\mu}$ be matrices, with $\mu=1,\ldots,d$,
and let $\bb v\in \CC^N$, with $N=n_1\cdots n_d$. Let 
$\bb V \in \CC^{n_1 \times \cdots \times n_d}$ be an order-$d$ tensor such
that $\mathrm{vec}(\bb V) = \bb v$, where $\mathrm{vec}$ denotes the operator
which stacks by columns the elements of the input tensor. Then, we have
\begin{equation*}
  (L_d \otimes L_{d-1} \otimes \cdots \otimes L_1)\bb v = 
  \mathrm{vec}(\bb V \times_1 L_1 \times_2 \cdots \times_d L_d).
\end{equation*}
\end{lemma}

As observed in the introduction for the case $d=2$,
since $K$ is a Kronecker sum of $d$ matrices,
we can similarly write
\begin{equation*}
  \exp(K)\bb v=
  \left(\exp(A_d)\otimes \exp(A_{d-1})\otimes\cdots\otimes \exp(A_1)\right)\bb v
\end{equation*}
which, using the introduced tensor formalism with $L_\mu=\exp(A_\mu)$,
can be computed as
\begin{equation}\label{eq:tucker}
  \exp(K)\bb v=\mathrm{vec}\left(\bb V \times_1 \exp(A_1)\times_2
  \cdots\times_d \exp(A_d)\right).
\end{equation}
The superiority of such an approach to compute
the action $\exp(K)\bb v$ has been
thoroughly analyzed and highlighted in Reference~\cite{CCEOZ22}
in the context of the numerical solution of some
Schr\"odinger equations.
Also, this technique has been successfully used in Reference~\cite{CCZ23kp}
for some (linear) advection--diffusion--absorption equations with
space dependent coefficients.

Formula~\eqref{eq:tucker} can be employed to compute
the exact solution of
system~\eqref{eq:ODE} in the case $\bb g(t,\bb u(t))\equiv 0$.
{\color{black}In general}, it is possible to integrate the
system} by employing an integrator of stiff order one
such as exponential Euler
\begin{subequations}\label{eq:expEul}
\begin{equation}\label{eq:expEuldp}
\bb u_{n+1}=\bb u_n+\tau\varphi_1(\tau K)\bb f(t_n,\bb u_n),
\end{equation}
where $\bb u_n\approx \bb u(t_n)$
and $\tau$ is the time step size, constant for simplicity of exposition.
Its implementation requires
the computation of the action of the $\varphi_1$ function on a vector.
An equivalent formulation of the scheme (see Reference~\cite{HO05}) is
\begin{equation}\label{eq:expEullc}
\bb u_{n+1}=\exp(\tau K)\bb u_n+\tau\varphi_1(\tau K)\bb g(t_n,\bb u_n),
\end{equation}
\end{subequations}
and a simple way to evaluate the latter is to compute the action
of the exponential of a slightly augmented matrix, that is
\begin{equation*}
  \exp\left(\begin{bmatrix}
    \tau K & \tau \bb g(t_n,\bb u_n)\\
    0\ \cdots\ 0 & 0
  \end{bmatrix}\right)\begin{bmatrix}
    \bb u_n\\
    1
    \end{bmatrix},
\end{equation*}
and take the first $N$ rows of the resulting vector.
The advantage of this approach is that it can be easily generalized to
higher order $\varphi$-functions~\cite{AMH11},
whose actions are needed for high stiff
order exponential integrators.
Indeed, consider the exponential
Runge--Kutta scheme of stiff order two
ETD2RK~\cite{CM02}
{\color{black}
\begin{subequations}\label{eq:ETD2RK}
  \begin{equation}\label{eq:ETD2RKdp}
  \begin{aligned}
    \bb u_{n2}&=\bb u_n+\tau\varphi_1(\tau K)\bb f(t_n,\bb u_n)\\
    \bb u_{n+1}&=\bb u_{n2}+
    \tau\varphi_2(\tau K)(\bb g(t_{n+1},\bb u_{n2})-\bb g(t_n,\bb u_n)),
    \end{aligned}
\end{equation}
or, equivalently,
  \begin{equation}\label{eq:ETD2RKlc}
  \begin{aligned}
    \bb u_{n2}&=\exp(\tau K)\bb u_n+\tau\varphi_1(\tau K)\bb g(t_n,\bb u_n)\\
    \bb u_{n+1}&=\exp(\tau K)\bb u_n+\tau\varphi_1(\tau K)\bb g(t_n,\bb u_n)+
    \tau\varphi_2(\tau K)(\bb g(t_{n+1},\bb u_{n2})-\bb g(t_n,\bb u_n)).
    \end{aligned}
  \end{equation}
  \end{subequations}}%
For its implementation, we can compute the actions of the two matrix functions
$\varphi_1(\tau K)$ and $\varphi_2(\tau K)$
{\color{black}(formulation~\eqref{eq:ETD2RKdp})}. Alternatively, we can compute
two linear combinations of actions of $\varphi$-functions
of type $\exp(\tau K)\bb v_0+\varphi_1(\tau K)\bb v_1+\varphi_2(\tau K)\bb v_2$
{\color{black}(formulation~\eqref{eq:ETD2RKlc})}.
To do the latter, it is again possible to use an
 augmented matrix approach.
 In fact, from Theorem~2.1 of Reference~\cite{AMH11},
we have that the first $N$ rows of the vector
{\color{black}
\begin{equation*}
  \exp\left(c\begin{bmatrix}
    &\tau K && \bb v_p & \bb v_{p-1} & \ldots & \bb v_2 & \bb v_1\\
    0 & \ldots & 0 & 0 & 1  & 0 & \ldots & 0\\[-1ex]
    \vdots & & \vdots & \vdots &  0 &  \ddots &  & \vdots \\[-1ex]
    \vdots &  & \vdots  & \vdots & \vdots  &  & \ddots & 0\\
    0 & \ldots & 0 & 0 & 0 & \ldots & 0 & 1\\
    0 & \ldots & 0 & 0 & 0 & \ldots & \ldots & 0
  \end{bmatrix}\right)\begin{bmatrix}\bb v_0\\
    0\\[-1ex]
    \vdots\\[-1ex]
        \vdots\\
    0\\
    1\end{bmatrix},\quad c\in\CC,
\end{equation*}
coincide with the vector
\begin{equation*}
  \exp(c\tau K)\bb v_0+c\varphi_1(c\tau K)\bb v_1+
  c^2\varphi_2(c\tau K)\bb v_2+\cdots+c^p\varphi_p(c\tau K)\bb v_p.
  \end{equation*}}%
With this technique, the computation of
a single linear combination of actions
of $\varphi$-functions
reduces to the action of the exponential of an augmented matrix.

However, in this context it is not possible to directly apply the
$\mu$-mode approach for computing the action
of the matrix exponential, since the augmented matrix is not anymore
a Kronecker sum.
The idea is then to take the integral definition
of the $\varphi$-functions~\eqref{eq:phifunctions}
and approximate it by a quadrature rule.
By doing so, for the exponential Euler scheme
introduced in formula~\eqref{eq:expEullc} we obtain
\begin{equation*}
  \bb u_{n+1}\approx \exp(\tau K)\bb u_n+\tau\sum_{i=1}^q w_i\exp((1-\theta_i)\tau K)
  \bb g(t_n,\bb u_n),
\end{equation*}
where $w_i$ and $\theta_i$ are the quadrature weights and nodes,
respectively (see Section~\ref{sec:sq} for an effective choice).

Since the matrix
arguments of the exponential function are Kronecker sums,
it is now possible to employ a $\mu$-mode approach for the efficient
evaluation of its action. In the next section,
we give more details and extend this idea to
higher order $\varphi$-functions and integrators.

\section[Approximation of actions of phi-functions
  of a Kronecker sum]{Approximation of $\varphi$-functions
  of a Kronecker sum}\label{sec:phiapprox}
In this section,
by using the tools presented in Section~\ref{sec:motivating},
we describe in detail how to approximate
the action of single \mbox{$\varphi$-functions} on the same vector and
the linear combination of actions of $\varphi$-functions,
{\color{black}which are the two tasks addressed by the algorithm \textsc{phiks}}.
For instance, a $\nu$-stage explicit exponential Runge--Kutta
integrator~\cite{HO10} with time step size $\tau$ is defined by
\begin{subequations}\label{eq:RK}
\begin{equation}
  \begin{aligned}
    \bb u_{ni}&=\exp(c_i\tau K)\bb u_n+c_i\tau \varphi_1(c_i\tau K)\bb g(t_n,\bb u_n)+
    \tau\sum_{j=2}^{i-1}a_{ij}(\tau K)\bb d_{nj}\\
    &=\bb u_n+c_i\tau \varphi_1(c_i\tau K)\bb f(t_n,\bb u_n)+
    \tau\sum_{j=2}^{i-1}a_{ij}(\tau K)\bb d_{nj},& 2&\le i\le \nu,\\
    \bb u_{n+1}&=\exp(\tau K)\bb u_n+\tau \varphi_1(\tau K)\bb g(t_n,\bb u_n)+
    \tau\sum_{i=2}^\nu b_i(\tau K)\bb d_{ni}\\
    &=\bb u_n+\tau \varphi_1(\tau K)\bb f(t_n,\bb u_n)+
    \tau\sum_{i=2}^\nu b_i(\tau K)\bb d_{ni},
  \end{aligned}
\end{equation}
where
\begin{equation}
  \bb d_{ni}=\bb g(t_n+c_i\tau,\bb u_{ni})-\bb g(t_n,\bb u_n).
\end{equation}
\end{subequations}
{\color{black}
Notice that for $\nu=1$ this integrator reduces to the exponential Euler
method~\eqref{eq:expEul}, while for $\nu=2$, $c_2=1$, and $b_2=\varphi_2$
we retrieve the ETD2RK method~\eqref{eq:ETD2RK}.
The generic scheme~\eqref{eq:RK} can be written in a compact way using the 
reduced tableau~\cite{LO14b} 
\begin{equation*}
  \begin{array}{c|cccc}
    c_2 &\\
    c_3 & a_{32}\\
    \vdots & \vdots & \ddots\\
    c_\nu & a_{\nu 2} & \cdots & a_{\nu \nu-1}\\
    \hline
     & b_2 & \cdots & b_{\nu-1} & b_\nu
    \end{array}
\end{equation*}
Here and throughout the paper,
by ``reduced tableau'' we mean that for each stage and for the final
approximation $\bb u_{n+1}$
we write only the coefficients
corresponding to the perturbation of the underlying exponential Euler scheme.}
The matrix functions $a_{ij}$ and $b_i$ are
linear
combinations of $\varphi$-functions.
Thus, each stage $\bb u_{ni}$ and the approximation $\bb u_{n+1}$
turn out to be {\color{black}\emph{general}} linear combinations of
actions of $\varphi$-functions, not necessarily
in form~\eqref{eq:lincomb}.
In fact,
it could be more convenient to compute actions of different
$\varphi$-functions on the same vector~\eqref{eq:diffphi}, and use them
to assemble the stages and the final approximation
{\color{black}(see the discussion in Section~\ref{sec:column-wise}).}
We start by describing this procedure.
\subsection[Actions of phi-functions on the same
  vector]{Actions of $\varphi$-functions on the same
  vector}\label{sec:diffphi}
Let us consider for the moment the simple second-order exponential
Runge--Kutta {\color{black}method~\eqref{eq:ETD2RKdp}, which
requires the computation of}
$\varphi_1(\tau K)\tau \bb f(t_n,\bb u_n)$
and $\varphi_2(\tau K)\tau (\bb g(t_{n+1},\bb u_{n2})-\bb g(t_n,\bb u_n))$.
We start with the computation of $\varphi_1(K)\bb v$, where, for clarity
of exposition, we omit the time step size $\tau$
and use a generic vector $\bb v$.
As mentioned in Section~\ref{sec:motivating}, the idea is to fully
exploit the possibility to apply the Tucker operator to compute
actions of suitable matrix exponentials.
Hence, we directly
approximate  the integral representation
\begin{equation}\label{eq:phi1quad}
\varphi_1(K)\bb v=\int_0^1 \exp((1-\theta)K)\bb vd\theta
\end{equation}
by a quadrature formula.
To avoid an impractical number of quadrature points,
we introduce
a scaling strategy. Therefore,
the quadrature rule is applied to the computation of
$\varphi_1(K/2^s)\bb v$,
that is
\begin{equation*}
  \varphi_1(K/2^s)\bb v=
  \int_0^1 \exp((1-\theta)K/2^s)\bb v d\theta\approx
  \sum_{i=1}^q w_i\exp((1-\theta_i)K/2^s)\bb v,
\end{equation*}
where $\theta_i$ and $w_i$ are $q$ quadrature nodes and weights,
respectively. Notice that
  we choose to scale the matrix $K$ by a power of two to
employ the favorable scaling and squaring algorithm~\cite{SW09} for
$\varphi$-functions. The choices of the quadrature formula, of
the number $q$ of quadrature nodes, and of the nonnegative
integer scaling $s$ will be discussed in detail in Section~\ref{sec:sq}.
Then, the evaluation  of
the integrand above at each quadrature point
$\theta_i\in[0,1]$ can be performed by the Tucker operator
\begin{equation}\label{eq:Tuckerprod}
    \bb V \times_1
    \exp((1-\theta_i)A_1/2^s)\times_2
    \cdots
  \times_d \exp((1-\theta_i) A_d/2^s),
  \end{equation}
see formula~\eqref{eq:tucker}.
Finally, to recover $\varphi_1(K)\bb v$ from its scaled version,
we use the following squaring formula (see again Reference~\cite{SW09})
\begin{equation*}
\left\{\begin{aligned}
  \varphi_1(K/2^{j-1}) \bb v&=\frac{1}{2}
    \left(\exp(K/2^j)\varphi_1(K/2^j) \bb v+
    \varphi_1(K/2^j) \bb v\right),\\
    \exp(A_\mu/2^{j-1})&=  \exp(A_\mu/2^j)\exp(A_\mu/2^j),
    \end{aligned}\right.
\end{equation*}
which has to be repeated for $j=s,s-1,\ldots,1$.
To perform the squaring,
\emph{no} full matrix $\exp(K/2^j)$ has to be evaluated
in practice. In fact, to compute its action on
$\varphi_1(K/2^j)\bb v$, which is available
as a tensor, it is enough to compute the Tucker operator
with the small sized matrices $\exp(A_\mu/2^j)$.
Notice that the idea of approximating
  the integral definition of the $\varphi$-functions
  by a quadrature formula and computing the action
  of the matrix exponential by a Tucker operator has also recently
  been presented in Reference~\cite{CMM23}.
  Nevertheless, we decided to report the above
  description for the sake of clarity and to introduce later
  additional features critical for exponential integrators, such as the
  effective usage of scaled quantities
  (see formula~\eqref{eq:scaldiffphi}) and the
  extension of the technique to linear combinations of actions of
  $\varphi$-functions (see the next section).
  
  Let us proceed by considering
  the approximation of the action of $\varphi_2(K)$, that
is
\begin{equation}\label{eq:phi2quad}
  \varphi_2(K)\bb v=\int_0^1 \theta\exp((1-\theta)K)
  \bb vd\theta.
\end{equation}
Comparing integrals \eqref{eq:phi1quad} and
\eqref{eq:phi2quad} it appears
clear that, if we define a common scaling strategy,
we can compute
the two approximations at once by selecting the same quadrature nodes
and weights, but different integrand functions
\begin{equation*}
  \exp((1-\theta)K/2^s)\bb v\quad \text{and}\quad
  \theta\exp((1-\theta)K/2^s)\bb v.
\end{equation*}
Therefore, the two quadrature formulas can be implemented
with common evaluations of
the matrices
$\exp((1-\theta_i)A_\mu/2^s)$
for each quadrature point $\theta_i$ and each $\mu$.
Their action on $\bb v$ is computed with a single
Tucker operator~\eqref{eq:Tuckerprod}, followed by
the multiplication by the scalar $\theta_i$ needed for
the approximation of $\varphi_2(K/2^s)\bb v$.
After assembling the quadrature,
the steps of the squaring are
\begin{equation*}
  \left\{
  \begin{aligned}
    \varphi_2(K/2^{j-1})\bb v&=
      \frac{1}{4}\left(\exp(K/2^j)\varphi_2(K/2^j)\bb v+\varphi_1(K/2^j) \bb v+
    \varphi_2(K/2^j)\bb v\right),\\
    \varphi_1(K/2^{j-1})\bb v&=\frac{1}{2}\left(\exp(K/2^j)
    \varphi_1(K/2^j)\bb v+
    \varphi_1(K/2^j)\bb v\right),\\
    \exp(A_\mu/2^{j-1})&=\exp(A_\mu/2^j)\exp(A_\mu/2^j),
\end{aligned}\right.
\end{equation*}%
to be repeated for $j=s,s-1,\ldots,1$.

The generalization to the
computation of the action of the first~$p$ $\varphi$-functions
on the \emph{same} vector $\bb v$
\begin{equation*}
  \{\varphi_1(K)\bb v,\varphi_2(K)\bb v,\ldots,  \varphi_p(K)\bb v\}
\end{equation*}
is straightforward. First, we compute their approximations at the same scaled
matrix by the common quadrature rule, \mbox{i.e.},
\begin{subequations}\label{eq:diffphicompl}
\begin{equation}\label{eq:integrandsimul}
  \varphi_\ell(K/2^s)\bb v\approx
    \sum_{i=1}^qw_i\frac{\theta_i^{\ell-1}}{(\ell-1)!}
    \exp((1-\theta_i)K/2^s)
  \bb v,\quad \ell=1,\ldots,p.
\end{equation}
Then, we perform the squaring procedure
\begin{equation}\label{eq:simul}
\left\{\begin{aligned}
  \varphi_\ell(K/2^{j-1})\bb v&=
    \frac{1}{2^\ell}\Big(\exp(K/2^j)\varphi_\ell(K/2^j)\bb v
    +
    \sum_{k=1}^\ell \frac{\varphi_k(K/2^j)\bb v}{(\ell-k)!}\Big),
    & \ell&=p,p-1,\ldots,1,\\
    \exp(A_\mu/2^{j-1})&=\exp(A_\mu/2^j)\exp(A_\mu/2^j),
  \end{aligned}\right.
\end{equation}%
\end{subequations}
for $j=s,s-1,\ldots,1$.
{\color{black}The action  of $\exp(K)$ on $\bb v$ can be computed using
  an additional single Tucker operator of type~\eqref{eq:tucker}.}
We stress that the relevant computations in
formulas~\eqref{eq:diffphicompl} are performed
by means of the $\mu$-mode approach,
without forming the large sized matrix $K$.
Notice also that, from {\color{black}formulas \eqref{eq:diffphicompl}, we
  obtain
at no additional cost also
$\varphi_\ell(K/2^{j-1})\bb v$, 
  $j=2,\ldots,\hat s$, where $\hat s\le s+1$ is the number of
  desired scales. This feature
could be useful for the efficient implementation of exponential
integrators that require, for instance, the quantities
\begin{equation}\label{eq:scaldiffphi}
  \exp(c_j K)\bb v,\ \varphi_1(c_j K)\bb v,\ 
  \varphi_2(c_j K)\bb v,\ \ldots,
  \varphi_p(c_j K)\bb v,\quad j=1,\ldots,\hat s,
  \end{equation}
with $c_j=c/2^{j-1}$ and $c\in\CC$},
as shown in the numerical examples of
  Sections~\ref{sec:ac} and \ref{sec:column-wise}.
\begin{remark}
Notice that the quadrature rule in formula~\eqref{eq:integrandsimul}
is equivalent to
\begin{equation*}
    \sum_{i=1}^qw_i\frac{\theta_i^{\ell-1}}{(\ell-1)!}
    \exp((1-\theta_i)(K-\sigma I)/2^s)
  \rme^{(1-\theta_i)\sigma/2^s}
  \bb v,\quad \ell=1,\ldots,p,
\end{equation*}
where $\sigma\in\CC$ is a shift parameter. Given the
Kronecker sum structure of $K$, it is possible
to choose $\sigma$ as
the sum of $d$ shifts $\sigma_\mu$, selected  in such
a way that $A_\mu-\sigma_\mu I$ has a smaller norm than $A_\mu$
(and thus its exponential can be possibly computed in a more efficient
way~\cite{AMH11,CZ19}).
A common and effective choice
for $\sigma_\mu$
is the trace of the matrix $A_\mu$ divided by $n_\mu$, which
corresponds to its average eigenvalue and minimizes the Frobenius
norm of $A_\mu-\sigma_\mu I$.
\end{remark}
We now summarize the number of Tucker
operators of the whole procedure {\color{black}inside \textsc{phiks}} needed
to obtain the quantities in formula~\eqref{eq:scaldiffphi}. We recall that,
if we assume $n_1=\ldots=n_d=n$, the
computational cost of a single Tucker operator
is $\mathcal{O}(n^{d+1})$. For each quadrature point
we need to compute one Tucker operator.
Then, for each step of the
squaring phase, we have $p$ Tucker operators.
Finally, we have one Tucker operator for the computation of
$\exp(K/2^{j-1})\bb v$
for each $j=1,\ldots,\hat s$.
Therefore, the total number of Tucker operators is
\begin{equation}\label{eq:Tpsv}
  T_\#(s,\hat s,q,p)=  q+sp+\hat s.
\end{equation}
We remark that for $d\ge 3$ the number $T_\#$ gives an
adequate indication of the computational cost of the whole procedure,
being the Tucker operator the most expensive operation.
On the other hand, for $d<3$ other tasks
such as the computation of the matrix exponential may have a comparable
cost (or even higher, in the trivial
case $d=1$).
\subsection[Linear combination of actions of
  phi-functions]{Linear combination of actions of
  $\varphi$-functions}\label{sec:lincomb}
Let us consider for the moment the
Runge--Kutta {\color{black}scheme~\eqref{eq:ETD2RKlc}, which
requires the two linear combinations
of actions of $\varphi$-functions
$\exp(\tau K)\bb u_n+\varphi_1(\tau K) \tau\bb g(t_n,\bb u_n)$ and
$\exp(\tau K)\bb u_n+\varphi_1(\tau K) \tau\bb g(t_n,\bb u_n)+
\varphi_2(\tau K) \tau(\bb g(t_{n+1},\bb u_{n2})-\bb g(t_n,\bb u_n))$.}

We then introduce the compact notation
 \begin{equation*}
  \Phi_j(K)\bb v_{i_1,i_2,\ldots,i_p} =
  \frac{\varphi_1(K/2^j)\bb v_{i_1}}{2^{j}}+
  \frac{\varphi_2(K/2^j)\bb v_{i_2}}{2^{2j}}+\cdots+
  \frac{\varphi_p(K/2^j)\bb v_{i_p}}{2^{pj}},\quad p>1,
\end{equation*}
with $j$ nonnegative integer
and, for simplicity of exposition, we describe in detail the approximation of
$\Phi_0(K)\bb v_{1,2}$.
The idea is to  apply a quadrature formula to the integral
\begin{equation*}
\Phi_0(K)\bb v_{1,2}=\int_0^1 \exp((1-\theta)K)(\bb v_1+\theta \bb v_2)d\theta
\end{equation*}
in combination with a scaling strategy. To do so,
we first approximate by
a common quadrature rule the scaled linear combinations
\begin{equation*}
  \begin{aligned}
    \Phi_s(K)\bb v_2&=\int_{0}^1\exp((1-\theta)K/2^s)\frac{\bb v_2}{2^s}d\theta
    \approx\sum_{i=1}^qw_i\exp((1-\theta_i)K/2^s)\frac{\bb v_2}{2^s},\\
    \Phi_s(K)\bb v_{1,2}&=\int_{0}^1\exp((1-\theta)K/2^s)\left(
    \frac{\bb v_1}{2^s}+\theta\frac{\bb v_2}{2^{2s}}\right)d\theta\\
    &\approx\sum_{i=1}^qw_i\exp((1-\theta_i)K/2^s)
\left(
    \frac{\bb v_1}{2^s}+\theta_i\frac{\bb v_2}{2^{2s}}\right).
  \end{aligned}
\end{equation*}
{\color{black}
Notice that the approximations in
the above formulas require the common evaluation of the matrices
$\exp((1-\theta_i)A_\mu/2^s)$
for every quadrature node $\theta_i$ and every $\mu$.
In addition, each of the two approximations needs
a single Tucker operator for every~$\theta_i$.
Then, it is possible to perform the squaring
procedure of $\Phi_s(K)\bb v_{1,2}$ by
\begin{equation*}
  \left\{\begin{aligned}
    \Phi_{j-1}(K)\bb v_{1,2}&=\exp(K/2^j)\Phi_j(K)\bb v_{1,2}
    +\frac{\Phi_j(K)\bb v_2}{2^j}+\Phi_j(K)\bb v_{1,2},\\
    \Phi_{j-1}(K)\bb v_2&=\exp(K/2^j)\Phi_j(K)\bb v_2+
    \Phi_j(K)\bb v_2,\\
    \exp(A_\mu/2^{j-1})&=\exp(A_\mu/2^j)\exp(A_\mu/2^j),
  \end{aligned}\right.
\end{equation*}%
with $j=s,s-1,\ldots,1$. The first two rows can be easily verified
by computing blocks $(1,3)$ and $(1,2)$ of the squares
of the matrices in equality
\begin{equation*}
    \exp\left(2^{-j}\begin{bmatrix}K&\bb v_2&\bb v_1\\
      0&0&1\\
    0&0&0\end{bmatrix}\right)=
  \begin{bmatrix}
    \exp(K/2^j) & \Phi_j(K) \bb v_2&
    \Phi_j(K) \bb v_{1,2}\\
    0 & 1 & 2^{-j}\\
    0 & 0 &1
  \end{bmatrix}.
\end{equation*}}%

The generalization to a \emph{linear combination} of
actions of the first~$p$ $\varphi$-functions on the vectors $\bb v_\ell$
\begin{equation*}
  \Phi_0(K)\bb v_{1,2,\ldots,p}=\varphi_1(K)\bb v_1+
  \varphi_2(K)\bb v_2+\cdots
  +\varphi_p(K)\bb v_p
\end{equation*}
requires first the application of the common quadrature rule
to the scaled linear combinations
\begin{subequations}\label{eq:lincombcompl}
\begin{equation}\label{eq:integrandcomb}
  \Phi_s(K)\bb v_{p-\ell+1,p-\ell+2,\ldots,p}\approx
  \sum_{i=1}^qw_i  \exp((1-\theta_i)K/2^s)\left(
  \sum_{k=1}^\ell\frac{\theta_i^{\ell-k}}{(\ell-k)!}
  \frac{\bb v_{p+1-k}}{2^{(\ell-k+1)s}}\right),
\end{equation}
with $\ell=1,\ldots,p$.
Then, the squaring procedure
\begin{equation}\label{eq:comb}
\left\{
  \begin{aligned}
    \Phi_{j-1}(K)\bb v_{p-\ell+1,p-\ell+2,\ldots,p}&=
    \exp(K/2^j)\Phi_j(K)\bb v_{p-\ell+1,p-\ell+2,\ldots,p}\\
    &+\sum_{k=1}^\ell\frac{\Phi_j(K)\bb v_{p-k+1,p-k+2,\ldots,p}}{(\ell-k)!2^{(\ell-k)j}},
    \quad \ell=p,p-1,\ldots,1,\\
        \exp(A_\mu/2^{j-1})&=\exp(A_\mu/2^j)\exp(A_\mu/2^j)
  \end{aligned}\right.
 \end{equation}%
\end{subequations}
has to be repeated for $j=s,s-1,\ldots,1$.
{\color{black}Clearly, the action of $\exp(K)$ on the
  vector $\bb v_0$ can be added to
the linear combination $\Phi_0(K)\bb v_{1,2,\ldots,p}$ at the cost
of a single additional Tucker operator {\color{black}of type \eqref{eq:tucker}}.}
We stress again that all the computations
in~\eqref{eq:lincombcompl} are performed
in a $\mu$-mode fashion.
Notice also that, from the squaring formula, we obtain at no additional
cost the
quantities $\Phi_{j-1}(K)\bb v_{1,2,\ldots,p}$, $j=2,\ldots, \hat s$,
which can be employed in the efficient implementation of 
exponential integrators that require, for instance, combinations
of the form
{\color{black}
\begin{equation}\label{eq:scallincomb}
  \exp(c_j K)\bb v_0+c_j\varphi_1(c_j K)\bb v_1+
  c_j^2\varphi_2(c_j K)\bb v_2+\cdots+
  c_j^p\varphi_p(c_j K)\bb v_p,\quad
  j=1,\ldots,\hat s,
\end{equation}%
with $c_j=c/2^{j-1}$ and $c\in\CC$.}
An explicit example will be presented in
the numerical experiment of Section~\ref{sec:row-wise}.

We now summarize the number of Tucker
operators needed by the whole procedure
{\color{black} inside \textsc{phiks}}
to obtain the quantities in formula~\eqref{eq:scallincomb}.
For each quadrature point
we need to compute $p$
Tucker operators. Then, for each step of the
squaring phase, we have $p$ Tucker operators.
Finally, we have one Tucker operator
for the computation of
$\exp(K/2^{j-1})\bb v_0$ for each $j=1,\ldots,\hat s$.
Therefore, the total number of Tucker operators is
\begin{equation}\label{eq:Tlcp}
  T_\#(s,\hat s,q,p)=
    qp+sp+\hat s.
\end{equation}
\subsection[Choice of s, q,
  and quadrature formula]{Choice of $s$, $q$,
  and quadrature formula}\label{sec:sq}
The choice of the scaling value $s$ and the number of
quadrature points $q$ is based on a suitable expansion of the error of the
quadrature formula.
After this selection,
the algorithm is \emph{direct} and no convergence test or
exit criterion is needed.
We start writing
\begin{equation}\label{eq:philrw}
  \varphi_\ell(K)=\int_{0}^{1}f_\ell(\theta,K)d\theta=
  \sum_{i=1}^{q}w_i f_\ell(\theta_i,K)+
                             R_{q}(f_\ell(\cdot,K)),
\end{equation}
where $R_{q}(f_\ell(\cdot,K))$ is the remainder
\begin{equation}\label{eq:remainder}
  R_{q}(f_\ell(\cdot,K)) =
  \frac{1}{2\pi\rmi}\oint_\Gamma k_{q}(z)f_\ell(z,K)dz,
\end{equation}
see Section~4.6 of Reference~\cite{DR84}.
Here, $\Gamma\subset\CC$ is an arbitrary simple closed curve surrounding
the interval $[0,1]$
and $k_{q}$ is the kernel defined by
\begin{equation*}
k_{q}(z)=\int_{0}^1\frac{\pi_q(t)}{\pi_q(z)(z-t)}dt,
\end{equation*}
with $\pi_q(t)$ the monic polynomial of degree $q$ with the quadrature
points as roots. {\color{black}Now, we describe a practical procedure
  to choose the scaling value $s$ and the number of quadrature points
  $q$ which balances accuracy and efficiency of the overall algorithms.
  In particular, we determine the parameters so that the remainder
  is below a certain tolerance, while the number of Tucker operators
  $T_\#$ is kept reasonably low to ensure that the computational
  cost is not excessive.}

Let us consider first the case of actions of $\varphi$-functions on the
same vector, as in formula~\eqref{eq:integrandsimul}. Then,
given a tolerance $\delta$ and starting from {\color{black}the
  scaling $s_0=0$},
  we look for the smallest number
  $q_0\in[q_\mathrm{min},q_\mathrm{max}]$ of
  quadrature points such that
\begin{equation*}
  \lVert R_{q_0}(f_\ell(\cdot,K))\rVert\rVert\bb v\rVert\le \delta,
  \quad \ell=1,\ldots,p.
\end{equation*}
We then repeat the procedure for increasing values of the
scaling $s_j\in\{1,2,\ldots\}$ and
look for the corresponding smallest value $q_j$ such that
\begin{equation*}
  \lVert R_{q_j}(f_\ell(\cdot,K/2^{s_j}))\rVert
  \lVert \bb v\rVert\le \delta \cdot 2^{\ell s_j},
  \quad \ell=1,\ldots,p.
\end{equation*}
Here, the tolerance is amplified by the factor $2^{\ell s_j}$ because
we take into account that squaring formula~\eqref{eq:simul}
requires $s_j$ divisions by $2^\ell$.
We continue until the
number of Tucker operators $T_\#(s_{\bar j+1},\hat s,q_{\bar j+1},p)$
in formula~\eqref{eq:Tpsv}
is larger than $T_\#(s_{\bar j},\hat s,q_{\bar j},p)$.
The obtained values $s=s_{\bar j}$ and $q=q_{\bar j}$
 are then employed in the
approximation of actions of $\varphi$-functions applied on
the vector $\bb v$ through formulas~\eqref{eq:diffphicompl}.

On the other hand, when considering a linear combination of actions
of $\varphi$-functions,
the quantities to be computed in
formula~\eqref{eq:integrandcomb}, starting from {\color{black}the
  scaling $s_0=0$},
correspond to the integrand functions 
\begin{equation*}
  \sum_{k=1}^\ell f_{\ell-k+1}(\cdot,K)\bb v_{p+1-k},\quad
  \ell=1,\ldots,p.
\end{equation*}
Therefore, for a given tolerance $\delta$,
we look for the smallest number
$q_0\in[q_\mathrm{min},q_\mathrm{max}]$
of
quadrature points
needed for all the values such that
\begin{equation*}
  \sum_{k=1}^\ell \lVert R_{q_0}(f_{\ell-k+1}(\cdot,K))\rVert
  \lVert\bb v_{p+1-k}\rVert\le \delta,\quad \ell=1,\ldots,p.
\end{equation*}
Then, similarly to the previous case,
we repeat the calculation for increasing values of the
scaling $s_j\in\{1,2,\ldots\}$
to obtain the corresponding smallest value $q_j$ such that
\begin{equation*}
  \sum_{k=1}^\ell \lVert R_{q_j}(f_{\ell-k+1}(\cdot,K/2^{s_j}))\rVert
  \frac{\lVert\bb v_{p+1-k}\rVert}{2^{(\ell-k+1)s_j}}\le \delta,
  \quad \ell=1,\ldots,p.
\end{equation*}
We continue this procedure
until the number of Tucker
operators $T_\#(s_{\bar j+1},\hat s,q_{\bar j+1},p)$ in formula~\eqref{eq:Tlcp}
is larger than $T_\#(s_{\bar j},\hat s,q_{\bar j},p)$.
The obtained values $s=s_{\bar j}$ and $q=q_{\bar j}$
are then employed in the
approximation of the linear combination of
$\varphi$-functions applied to the vectors $\bb v_1$, \ldots,
$\bb v_p$ through formulas~\eqref{eq:lincombcompl}.

The previous estimates clearly require computable
bounds for the remainders with different numbers of quadrature points,
integrand functions and scaling parameters.
To avoid cumbersome notation, we explain
the procedure for $R_q(f_\ell(\cdot,K))$ in formula~\eqref{eq:remainder}.
We choose $\Gamma=\Gamma_r$ to be the ellipse with foci in
$\{0,1\}$ and logarithmic capacity (half sum of its semi-axes) $r>1/4$,
that is
\begin{equation*}
  \Gamma_r=\left\{z\in\CC:z=z(\zeta)=r\rme^{\rmi\zeta}+\frac{1}{2}+
  \frac{\rme^{-\rmi\zeta}}{16r},\ \text{with $\zeta\in[0,2\pi)$} \right\}.
\end{equation*}
Then, we have
\begin{equation*}
  \begin{split}
    \lVert R_{q}(f_\ell(\cdot,K))\Vert &=
\Bigg\lVert \frac{1}{2\pi\rmi}\oint_{\Gamma_r} k_{q}(z)f_\ell(z,K)dz \Bigg\rVert
   =\\
   &=\frac{1}{2\pi}\Bigg\lVert \int_0^{2\pi}k_{q}(z(\zeta))
   f_\ell(z(\zeta),K)
   \left(r\rme^{\rmi \zeta}-
        \frac{\rme^{-\rmi \zeta}}{16r}\right)d\zeta \Bigg\rVert.
\end{split}
  \end{equation*}%
Finally, by using the fact that the numerical range of $K$
(denoted by $\mathcal{W}(K)$),
is a
$(1+\sqrt{2})$-spectral
set~\cite{CP17}, we estimate in 2-norm
\begin{equation}\label{eq:CPestimate}
  \lVert R_{q}(f_\ell(\cdot,K))\Vert_2 \leq
  \frac{1+\sqrt{2}}{2\pi}\sup_{w\in\Omega}\Bigg\lvert
\int_0^{2\pi}k_{q}(z(\zeta))
   f_\ell(z(\zeta),w)
   \left(r\rme^{\rmi \zeta}-
        \frac{\rme^{-\rmi \zeta}}{16r}\right)d\zeta
  \Bigg\rvert,
\end{equation}%
being $\Omega\subset\CC$ a smooth, bounded, convex domain which
embraces $\mathcal{W}(K$).
In our situation, we can easily find such a domain
without assembling the matrix $K$.
Indeed, it is possible to show that
\begin{equation*}
  \mathcal{W}(K)= \mathcal{W}(A_{\otimes 1})+\mathcal{W}(A_{\otimes 2})+\cdots
  +\mathcal{W}(A_{\otimes d})=
  \mathcal{W}(A_{1})+\mathcal{W}(A_{2})+\cdots+\mathcal{W}(A_{d}),
\end{equation*}
and $\mathcal{W}(A_\mu)$ can be estimated~\cite{CCZ20} with a rectangle
$\Xi_\mu$
obtained by computing the norms of
the Hermitian and the skew-Hermitian parts of the small sized matrices
$A_\mu$.
Thus, the rectangle $\Xi=\Xi_1+\cdots+\Xi_d$ embraces
$\mathcal{W}(K)$ and,
thanks to the maximum modulus principle,
the supremum in estimate~\eqref{eq:CPestimate} is attained at
the boundary of $\Xi$, which we suitably discretize.
Moreover, we approximate the integral by the trapezoidal rule.

Concerning the choice of the main quadrature formula~\eqref{eq:philrw},
we use the Gauss--Lobatto--Legendre one
{\color{black} with a number of quadrature points in the interval
  $[q_\mathrm{min},q_\mathrm{max}]=[3,12]$. These bounds
appear to be adequate for the addressed numerical experiments.}
Besides being
very accurate, it employs the endpoints of the integration interval $[0,1]$.
This allows on one side to
avoid one Tucker operator of type~\eqref{eq:Tuckerprod}
(since $\theta_q=1$), and on the other
to avail of the quantities $\exp(A_\mu/2^s)$
(corresponding to $\theta_1=0$), which
are needed for the squaring procedures.
{\color{black}This shrewdness, together with the fact that we can avoid computing
the Tucker operators for the matrix exponential in formula~\eqref{eq:scallincomb} if
$\bb v_0$ is zero, is taken into account in the actual implementation of
the algorithm {\color{black}\textsc{phiks}}.}
Finally,
the evaluation of the kernel $k_{q}$ in estimate~\eqref{eq:CPestimate}
is obtained by the recurrence relation
of the underlying orthogonal polynomials (see Reference~\cite{GV83}).
\section{Numerical experiments}\label{sec:numerical}
In this section, we validate our
MATLAB\footnote{The code is available at \url{https://github.com/caliarim/phiks}
  and is fully compatible with
  GNU Octave} implementation of \textsc{phiks} and
present the effectiveness of the proposed algorithm for the numerical
solution of stiff systems of ODEs with exponential Runge--Kutta
integrators from stiff order one to four. The implemented algorithm,
which works in \emph{any} space dimension $d$,
employs the function \textsc{tucker}
(contained in the package
KronPACK\footnote{\url{https://github.com/caliarim/KronPACK}})
to compute the underlying Tucker operators by means of $\mu$-mode products.
In addition, it
uses the internal MATLAB
function \textsc{expm} for the approximation
of the needed matrix exponentials.
Such a function is based on the double precision scaling and squaring
Pad\'e algorithm~\cite{AMH09}.

Concerning the two-dimensional example described in Section~\ref{sec:ac},
we compare the efficiency of our approach with
a technique recently introduced for the computation
of $\varphi$-functions of matrices that have Kronecker sum
structure~\cite{MMPC22}.
The method,
which was developed for the two-dimensional case only,
is direct, does not require an input tolerance, and
retrieves the action of a single $\varphi$-function of order $\ell$
by solving $\ell$ Sylvester equations.
This approach has some restrictions on the input matrices
($A_1$ and $-A_2$ must have disjoint spectra 
to have a unique solution of the Sylvester equation) and
it may suffer
of ill-conditioning for $\varphi$-functions of high order.
{\color{black}The accompanying
  software\footnote{\url{https://github.com/jmunoz022/Kronecker_EI}}
  of Reference~\cite{MMPC22} contains the scripts used by the
  authors to perform their
    numerical examples.
    For our purposes, we extracted the parts devoted to the
    computation of the $\varphi$-functions using the MATLAB
    function \textsc{sylvester}, and collected them in
    a function (that we named \textsc{sylvphi}) to be easily employed
    in our numerical experiments.}
In addition, in this two-dimensional example and in all the remaining
three-dimensional ones, we compare our approach with
recent and popular algorithms for
computing linear combinations of actions of $\varphi$-functions
for large and sparse general matrices, whose code
  in publicly available.
For convenience of the reader, we briefly describe them in the following.
\begin{itemize}
\item \textsc{phipm\_simul\_iom}\footnote{\url{https://github.com/drreynolds/Phipm_simul_iom}}
is a
Krylov subspace  solver with incomplete orthogonalization~\cite{LPR19}
which computes linear combinations of actions of $\varphi$-functions
at different time scales, by expressing everything in terms of the highest
order $\varphi$-function
and using a recurrence relation.
\item \textsc{kiops}\footnote{\url{https://gitlab.com/stephane.gaudreault/kiops/-/tree/master/}}
is another adaptive Krylov subspace solver
with incomplete orthogonalization~\cite{GRT18}.
It computes linear combinations
  of actions of $\varphi$-functions at different time scales
  by using the augmented matrix technique.
\item \textsc{bamphi}\footnote{\url{https://github.com/francozivcovich/bamphi}}
  is a hybrid Krylov-polynomial
  method~\cite{CCZ23}
  for computing 
  linear combinations of actions of $\varphi$-functions
  at different time scales,
  equipped with a backward error analysis of the underlying
  polynomial approximation. In contrast to the previous
  methods, it does not require to store a Krylov subspace.
\end{itemize}
We used all these methods with an incomplete orthogonalization
procedure of length two. Moreover, since
their MATLAB implementations output some information that can be
effectively used for successive calls, such as an estimate of the
appropriate Krylov subspace size, in our numerical experience we
obtained overall the best results by adopting the following
strategy: for each call of the routine at a certain time step we input the
information obtained by the same call at the previous time step.
In addition, these three methods, together with \textsc{phiks}, require
an input tolerance,
but their error estimates are substantially different.
For this reason, we decided to set the tolerance of each method to a
value proportional
to both the local error of
the used time marching scheme and the 2-norm of the current approximation
$\bb u_n$.
The proportionality constant has been selected for each method and each
integrator
as large as possible among the powers of two,
in such a way that the final error measured with respect to a reference
solution is not affected by the approximation error of the matrix functions.
We believe that running the experiments with tolerances obtained
in this way yields a fair comparison among all the
methods, ensuring the minimal effort needed to reach the accuracy of the
considered time marching schemes.
The study of a more sophisticated
technique for an effective choice of the tolerances
is far beyond the scope of this manuscript.

All the numerical experiments were performed on an
Intel\textsuperscript{\textregistered}
Core\textsuperscript{\texttrademark} i7-10750H
CPU with six physical cores and 16GB of RAM, using
MathWorks MATLAB\textsuperscript{\textregistered}
R2022a. The errors were measured in infinity norm relatively
to either the analytical solution, when available, or
to a reference solution computed by the fourth-order
integrator~\eqref{eq:expRK4s6}, implemented with the \textsc{phiks} routine
and a sufficiently large number of time steps.
 \subsection{Code validation}
 We extensively tested the \textsc{phiks} routine
 and we present here the results regarding the approximation of
 actions of $\varphi$-functions on the same vector and linear combinations
 of actions of $\varphi$-functions up to order $p=5$. The test
 matrices
 arise from the discretization by standard second
order finite differences of the complex operator
$(1+\rmi)/100\cdot\Delta$ in the spatial domain
$[0,1]^d$, for $d=3$ and $d=6$, with homogeneous Dirichlet boundary conditions.
The application
vectors $\bb v_1=\ldots=\bb v_p=\bb v$ are the
discretization of
\begin{equation*}
4096(1+\rmi) \prod_{\mu=1}^dx_\mu(1-x_\mu).
\end{equation*}
The number $n$ of discretization points for each spatial direction ranges from
64 to 121 for $d=3$, and from 8 to 11 for $d=6$.
As a term of comparison we consider the results obtained with
\textsc{kiops}. Both routines were called with input tolerance
set to the double precision unit roundoff value $2^{-53}$.
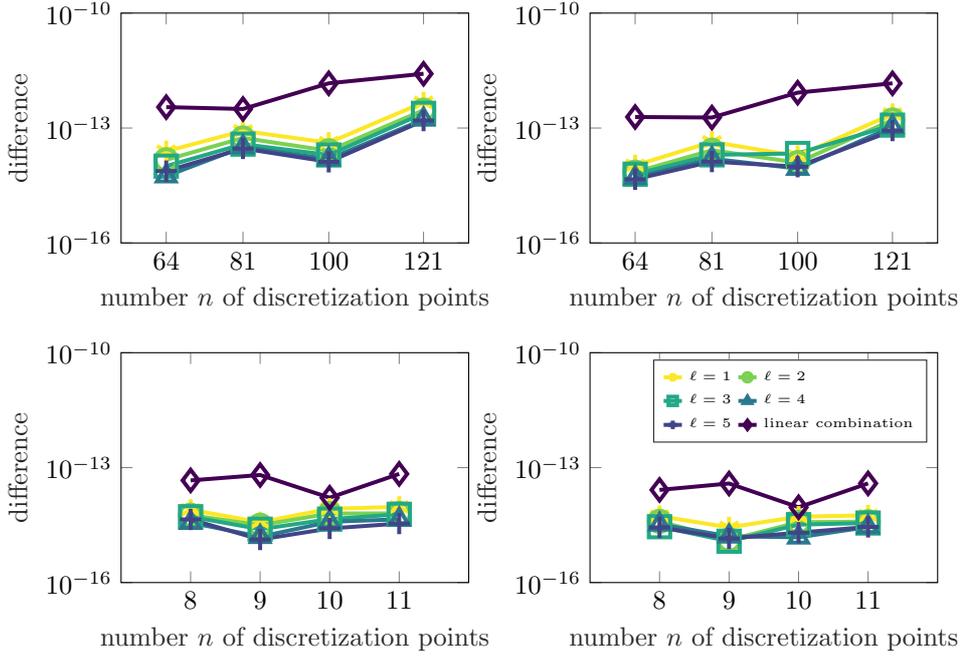
\begin{figure}[!ht]
  \centering
%
%
\definecolor{mycolor1}{rgb}{0.9922,0.9059,0.1451}%
\definecolor{mycolor2}{rgb}{0.4784,0.8196,0.3176}%
\definecolor{mycolor3}{rgb}{0.13330,0.65880,0.51760}%
\definecolor{mycolor4}{rgb}{0.16470,0.47060,0.55690}%
\definecolor{mycolor5}{rgb}{0.25490,0.26670,0.52940}%
\definecolor{mycolor6}{rgb}{0.26670,0.00400,0.32940}%
\begin{tikzpicture}

\begin{axis}[%
width=1.8in,
height=1.2in,
at={(0.822in,2.56in)},
scale only axis,
xmin=54,
xmax=131,
xlabel style={font=\color{white!15!black}},
xlabel={number $n$ of discretization points},
xtick={64,81,100,121},
xticklabels={{64},{81},{100},{121}},
ymode=log,
ymin=1e-16,
ymax=1e-10,
yminorticks=true,
ylabel style={font=\color{white!15!black}},
ylabel={difference},
axis background/.style={fill=white},
title style={font=\bfseries},
]
\addplot [color=mycolor1, line width=1.5pt, mark size=4pt, mark=asterisk, mark options={solid, mycolor1}]
  table[row sep=crcr]{%
64	2.47e-14\\
81      8.42e-14\\
100	4.13e-14\\
121	4.66e-13\\
};

\addplot [color=mycolor2, line width=1.5pt, mark size=4pt, mark=o, mark options={solid, mycolor2}]
  table[row sep=crcr]{%
64	1.48e-14\\
81	5.51e-14\\
100	2.55e-14\\
121	2.92e-13\\
};

\addplot [color=mycolor3, line width=1.5pt, mark size=4pt, mark=square, mark options={solid, mycolor3}]
  table[row sep=crcr]{%
64	1.01e-14\\
81	3.79e-14\\
100	1.98e-14\\
121	2.42e-13\\
};

\addplot [color=mycolor4, line width=1.5pt, mark size=4pt, mark=triangle, mark options={solid, mycolor4}]
  table[row sep=crcr]{%
64	5.56e-15\\
81	3.18e-14\\
100	1.54e-14\\
121	1.68e-13\\
};

\addplot [color=mycolor5, line width=1.5pt, mark size=4pt, mark=+, mark options={solid, mycolor5}]
  table[row sep=crcr]{%
64	7.53e-15\\
81	2.92e-14\\
100	1.31e-14\\
121	1.56e-13\\
};

\addplot [color=mycolor6, line width=1.5pt, mark size=4pt, mark=diamond, mark options={solid, mycolor6}]
  table[row sep=crcr]{%
64	3.51e-13\\
81	3.14e-13\\
100	1.47e-12\\
121	2.61e-12\\
};

\end{axis}

\begin{axis}[%
width=1.8in,
height=1.2in,
at={(3.25in,2.56in)},
scale only axis,
xmin=54,
xmax=131,
xlabel style={font=\color{white!15!black}},
xlabel={number $n$ of discretization points},
xtick={64,81,100,121},
xticklabels={{64},{81},{100},{121}},
ymode=log,
ymin=1e-16,
ymax=1e-10,
yminorticks=true,
ylabel style={font=\color{white!15!black}},
ylabel={difference},
axis background/.style={fill=white},
title style={font=\bfseries},
legend style={legend cell align=left, align=left, draw=white!15!black}
]
\addplot [color=mycolor1, line width=1.5pt, mark size=4pt, mark=asterisk, mark options={solid, mycolor1}]
  table[row sep=crcr]{%
64	1.06e-14\\
81	4.35e-14\\
100	1.85e-14\\
121	2.32e-13\\
};

\addplot [color=mycolor2, line width=1.5pt, mark size=4pt, mark=o, mark options={solid, mycolor2}]
  table[row sep=crcr]{%
64	7.06e-15\\
81	2.63e-14\\
100	1.25e-14\\
121	1.62e-13\\
};

\addplot [color=mycolor3, line width=1.5pt, mark size=4pt, mark=square, mark options={solid, mycolor3}]
  table[row sep=crcr]{%
64	6.22e-15\\
81	1.99e-14\\
100	2.17e-14\\
121	1.19e-13\\
};

\addplot [color=mycolor4, line width=1.5pt, mark size=4pt, mark=triangle, mark options={solid, mycolor4}]
  table[row sep=crcr]{%
64	5.05e-15\\
81	1.58e-14\\
100	8.75e-15\\
121	1.07e-13\\
};

\addplot [color=mycolor5, line width=1.5pt, mark size=4pt, mark=+, mark options={solid, mycolor5}]
  table[row sep=crcr]{%
64	4.59e-15\\
81	1.33e-14\\
100	9.79e-15\\
121	8.56e-14\\
};

\addplot [color=mycolor6, line width=1.5pt, mark size=4pt, mark=diamond, mark options={solid, mycolor6}]
  table[row sep=crcr]{%
64	1.94e-13\\
81	1.88e-13\\
100	8.38e-13\\
121	1.47e-12\\
};

\end{axis}

\begin{axis}[%
width=1.8in,
height=1.2in,
at={(0.822in,0.787in)},
scale only axis,
xmin=7,
xmax=12,
xlabel style={font=\color{white!15!black}},
xlabel={number $n$ of discretization points},
xtick={8,9,10,11},
xticklabels={{8},{9},{10},{11}},
ymode=log,
ymin=1e-16,
ymax=1e-10,
yminorticks=true,
ylabel style={font=\color{white!15!black}},
ylabel={difference},
axis background/.style={fill=white},
title style={font=\bfseries},
legend style={legend cell align=left, align=left, draw=white!15!black}
]
\addplot [color=mycolor1, line width=1.5pt, mark size=4pt, mark=asterisk, mark options={solid, mycolor1}]
  table[row sep=crcr]{%
8	7.95e-15\\
9	3.75e-15\\
10	8.45e-15\\
11	9.46e-15\\
};

\addplot [color=mycolor2, line width=1.5pt, mark size=4pt, mark=o, mark options={solid, mycolor2}]
  table[row sep=crcr]{%
8	5.75e-15\\
9	3.13e-15\\
10	6.33e-15\\
11	6.39e-15\\
};

\addplot [color=mycolor3, line width=1.5pt, mark size=4pt, mark=square, mark options={solid, mycolor3}]
  table[row sep=crcr]{%
8	5.24e-15\\
9	2.47e-15\\
10	4.51e-15\\
11	6.06e-15\\
};

\addplot [color=mycolor4, line width=1.5pt, mark size=4pt, mark=triangle, mark options={solid, mycolor4}]
  table[row sep=crcr]{%
8	3.66e-15\\
9	1.57e-15\\
10	3.85e-15\\
11	4.47e-15\\
};

\addplot [color=mycolor5, line width=1.5pt, mark size=4pt, mark=+, mark options={solid, mycolor5}]
  table[row sep=crcr]{%
8 4.44e-15\\
9 1.33e-15\\
10 2.58e-15\\
11 3.47e-15\\
};


\addplot [color=mycolor6, line width=1.5pt, mark size=4pt, mark=diamond, mark options={solid, mycolor6}]
  table[row sep=crcr]{%
8	4.64e-14\\
9	6.42e-14\\
10	1.64e-14\\
11	6.86e-14\\
};

\end{axis}

\begin{axis}[%
width=1.8in,
height=1.2in,
at={(3.25in,0.787in)},
scale only axis,
xmin=7,
xmax=12,
xlabel style={font=\color{white!15!black}},
xlabel={number $n$ of discretization points},
xtick={8,9,10,11},
xticklabels={{8},{9},{10},{11}},
ymode=log,
ymin=1e-16,
ymax=1e-10,
yminorticks=true,
ylabel style={font=\color{white!15!black}},
ylabel={difference},
axis background/.style={fill=white},
title style={font=\bfseries},
legend pos=north east,
legend style={legend cell align=left, align=left, draw=white!15!black, font=\tiny,
legend columns=2},
legend image post style={scale=0.4}
]
\addplot [color=mycolor1, line width=1.5pt, mark size=4pt, mark=asterisk, mark options={solid, mycolor1}]
  table[row sep=crcr]{%
8	5.41e-15\\
9	2.75e-15\\
10	5.25e-15\\
11	5.67e-15\\
};
\addlegendentry{$\ell = 1$}

\addplot [color=mycolor2, line width=1.5pt, mark size=4pt, mark=o, mark options={solid, mycolor2}]
  table[row sep=crcr]{%
8	3.96e-15\\
9	1.29e-15\\
10	3.76e-15\\
11	3.88e-15\\
};
\addlegendentry{$\ell = 2$}

\addplot [color=mycolor3, line width=1.5pt, mark size=4pt, mark=square, mark options={solid, mycolor3}]
  table[row sep=crcr]{%
8	2.87e-15\\
9	1.20e-15\\
10	3.25e-15\\
11	3.58e-15\\
};
\addlegendentry{$\ell = 3$}

\addplot [color=mycolor4, line width=1.5pt, mark size=4pt, mark=triangle, mark options={solid, mycolor4}]
  table[row sep=crcr]{%
8	3.31e-15\\
9	1.61e-15\\
10	1.52e-15\\
11	2.82e-15\\
};
\addlegendentry{$\ell = 4$}

\addplot [color=mycolor5, line width=1.5pt, mark size=4pt, mark=+, mark options={solid, mycolor5}]
  table[row sep=crcr]{%
8	2.75e-15\\
9	1.41e-15\\
10	2.01e-15\\
11	2.83e-15\\
};
\addlegendentry{$\ell = 5$}

\addplot [color=mycolor6, line width=1.5pt, mark size=4pt, mark=diamond, mark options={solid, mycolor6}]
  table[row sep=crcr]{%
8	2.62e-14\\
9	3.85e-14\\
10	9.24e-15\\
11	3.88e-14\\
};
\addlegendentry{linear combination}

\end{axis}

\end{tikzpicture}%
  \caption{Relative difference between
    \textsc{kiops} and \textsc{phiks}, measured in infinity norm,  for the
    actions $\varphi_\ell(K/2^j)\bb v$, $\ell=1,\ldots,p$,
    and for the linear combinations
    $\Phi_j(K)\bb v_{1,2,\ldots,p}$ in the code validation.
    The plots refer to $j=0$ and $d=3$ (top left),
    $j=1$ and $d=3$ (top right), $j=0$ and $d=6$ (bottom left),
    $j=1$ and $d=6$ (bottom right).}
    \label{fig:validation}
\end{figure}
We report in Figure~\ref{fig:validation} 
the relative difference in infinity norm between the approaches
and, for \textsc{phiks}, we collect in Table~\ref{tab:validation}
the values of the scaling parameter
$s$, the number of quadrature points $q$, and
the number of Tucker operators
{\color{black}$T_\#$ (see formulas~\eqref{eq:Tpsv} and \eqref{eq:Tlcp})}.
Overall, we observe an homogeneous behavior of the relative difference
between \textsc{kiops} and \textsc{phiks}
for all the values of $d$, $n$, and $\ell$, and a number of
Tucker operators required by the routine \textsc{phiks} which
increases very slowly with $n$.
\begin{table}[!ht]
  \centering
  \caption{Values of the scaling parameter {\color{black}$s$}, 
           number of quadrature points {\color{black}$q$},
    and number of Tucker operators {\color{black}$T_\#$} to compute actions of $\varphi$-functions
    on the same vector
    (top) and linear
    combinations of actions of $\varphi$-functions (bottom),
    employed by \textsc{phiks} in the code validation.}
\label{tab:validation}
  \begin{tabular}{c|cccc|cccc}
\multicolumn{9}{c}{$\varphi$-functions on the same vector }\\
\hline
& \multicolumn{4}{c|}{$d=3$} & 
    \multicolumn{4}{c}{$d=6$}\\
    \hline
$n$ & 64 & 81 & 100 & 121 & 8 & 9 & 10 & 11\\
\hline
    $s$ & 8 & 8 & 9 & 9 & 3 & 3 & 3 & 4\\
    $q$ & 10 & 12 & 11 & 12 & 11 & 11 & 12 & 10\\
    {\color{black}$T_\#$} & {\color{black}52} & {\color{black}54} & {\color{black}58} & 
    {\color{black}59} & {\color{black}28} & {\color{black}28} & {\color{black}29} & {\color{black}32}\\
  \end{tabular}
  \begin{tabular}{c|cccc|cccc}
\multicolumn{9}{c}{linear combination of $\varphi$-functions}\\
\hline
& \multicolumn{4}{c|}{$d=3$} & 
    \multicolumn{4}{c}{$d=6$}\\
    \hline
$n$ & 64 & 81 & 100 & 121 & 8 & 9 & 10 & 11\\
    \hline
    $s$ & 10 & 11 & 11 & 12 & 5 & 6 & 6 & 6\\
    $q$ & 7 & 7 & 8 & 7 & 8 & 7 & 7 & 7\\
    {\color{black}$T_\#$} & {\color{black}87} & {\color{black}92} & {\color{black}97} & 
    {\color{black}97} & {\color{black}67} & {\color{black}67} & {\color{black}67} & {\color{black}67}\\
  \end{tabular}
\end{table}

\subsection{Evolutionary
  advection--diffusion--reaction equation}\label{sec:adr}
In this section we consider the following evolutionary
advection--diffusion--reaction (ADR) equation 
\begin{subequations}\label{eq:ADR}
  \begin{equation}
    \left\{\begin{aligned}
  \partial_t u(t,x_1,x_2,x_3)&=\varepsilon \Delta u(t,x_1,x_2,x_3)+
  \alpha(\partial_{x_1} +\partial_{x_2}+\partial_{x_3})u(t,x_1,x_2,x_3)\\
  &+g(t,x_1,x_2,x_3,u(t,x_1,x_2,x_3)),\\
  u_0(x_1,x_2,x_3) &= 64x_1(1-x_1)x_2(1-x_2)x_3(1-x_3),
\end{aligned}\right.
    \end{equation}
in the spatial domain $[0,1]^3$, where the nonlinear function $g$
is defined by
\begin{equation}
  g(t,x_1,x_2,x_3,u(t,x_1,x_2,x_3))=\frac{1}{1+u(t,x_1,x_2,x_3)^2}+
  \Psi(t,x_1,x_2,x_3).
\end{equation}
Here, $\Psi(t,x_1,x_2,x_3)$ is chosen so that the analytical
solution is $u(t,x_1,x_2,x_3)=\rme^t u_0(x_1,x_2,x_3)$.
\end{subequations}
Finally, the equation is coupled with homogeneous Dirichlet boundary conditions.
The diffusion and advection
parameters are set to $\varepsilon=0.5$ and $\alpha=10$, respectively.
\begin{figure}[htb!]
  \centering
%
%
\definecolor{mycolor1}{rgb}{0.9922,0.9059,0.1451}%
\definecolor{mycolor2}{rgb}{0.3686,0.7882,0.38430}%
\definecolor{mycolor3}{rgb}{0.12940,0.56860,0.54900}%
\definecolor{mycolor4}{rgb}{0.23140,0.32160,0.54510}%
\definecolor{mycolor5}{rgb}{0.26670,0.00390,0.32940}%
\begin{tikzpicture}

\begin{axis}[%
width=1.8in,
height=1.5in,
at={(0.758in,0.481in)},
scale only axis,
xmode=log,
xmin=250,
xmax=840,
xminorticks=true,
xlabel style={font=\color{white!15!black}},
xlabel={number of time steps},
xtick={300,400,500,600,700},
xticklabels={{300},{400},{500},{600},{700}},
ymode=log,
ymax=1e-3,
ymin=1e-5,
yminorticks=true,
ylabel style={font=\color{white!15!black}},
ylabel={error},
axis background/.style={fill=white},
title style={font=\bfseries},
legend pos=south west,
legend style={legend cell align=left, align=left, draw=white!15!black, font=\tiny},
legend image post style={scale=0.4}
]
\addplot [color=mycolor2, line width=1.5pt, only marks, mark size=4pt, mark=asterisk, mark options={solid, mycolor2}]
  table[row sep=crcr]{%
300	0.000165235186389795\\
400	0.000123898171318099\\
500	9.91049913730389e-05\\
600	8.25799671443656e-05\\
700	7.07782213075619e-05\\
};
\addlegendentry{\textsc{kiops}}

\addplot [color=mycolor3, line width=1.5pt, only marks, mark size=4pt, mark=+, mark options={solid, mycolor3}]
  table[row sep=crcr]{%
300	0.000165235196540974\\
400	0.000123898173293434\\
500	9.91049919155734e-05\\
600	8.25799673310767e-05\\
700	7.07782213797785e-05\\
};
\addlegendentry{\textsc{bamphi}}

\addplot [color=mycolor4, line width=1.5pt, only marks, mark size=4pt, mark=o, mark options={solid, mycolor4}]
  table[row sep=crcr]{%
300	0.000165235244516632\\
400	0.000123898203193513\\
500	9.91050099112496e-05\\
600	8.25799787168147e-05\\
700	7.07782289689893e-05\\
};
\addlegendentry{\textsc{phipm\_simul\_iom}}

\addplot [color=mycolor5, line width=1.5pt, only marks, mark size=4pt, mark=square, mark options={solid, mycolor5}]
  table[row sep=crcr]{%
300	0.000165235138378534\\
400	0.000123898156126731\\
500	9.91049851245879e-05\\
600	8.25799641644721e-05\\
700	7.0778219642131e-05\\
};
\addlegendentry{\textsc{phiks}}

\addplot [color=black, dashed, line width=1.5pt]
  table[row sep=crcr]{%
300	0.000165149183219483\\
700	7.07782213797785e-05\\
};
\addlegendentry{order one}

\end{axis}

\end{tikzpicture}%
\hfill
%
%
\definecolor{mycolor1}{rgb}{0.9922,0.9059,0.1451}%
\definecolor{mycolor2}{rgb}{0.3686,0.7882,0.38430}%
\definecolor{mycolor3}{rgb}{0.12940,0.56860,0.54900}%
\definecolor{mycolor4}{rgb}{0.23140,0.32160,0.54510}%
\definecolor{mycolor5}{rgb}{0.26670,0.00390,0.32940}%
\begin{tikzpicture}

\begin{axis}[%
width=1.8in,
height=1.5in,
at={(0.758in,0.481in)},
scale only axis,
xmin=54,
xmax=131,
xminorticks=true,
xlabel style={font=\color{white!15!black}},
xlabel={number $n$ of discretization points},
xtick={64,81,100,121},
xticklabels={{64},{81},{100},{121}},
ymode=log,
ymin=1,
ymax=1000,
yminorticks=true,
ylabel style={font=\color{white!15!black}},
ylabel={wall-clock time},
axis background/.style={fill=white},
title style={font=\bfseries},
legend style={legend cell align=left, align=left, draw=white!15!black, font=\footnotesize, at={(0.54,0.98)}},
legend image post style={scale=0.6}
]
\addplot [color=mycolor5, line width=1.5pt, mark size=4pt, mark=square, mark options={solid, mycolor5}]
  table[row sep=crcr]{%
64	2.656309\\
81	6.35685\\
100	16.815645\\
121	34.748827\\
};

\addplot [color=mycolor3, line width=1.5pt, mark size=4pt, mark=+, mark options={solid, mycolor3}]
  table[row sep=crcr]{%
64	10.896792\\
81	27.403405\\
100	62.183351\\
121	140.971347\\
};

\addplot [color=mycolor2, line width=1.5pt, mark size=4pt, mark=asterisk, mark options={solid, mycolor2}]
  table[row sep=crcr]{%
64	19.773262\\
81	53.804316\\
100	133.122662\\
121	320.582258\\
};

\addplot [color=mycolor4, line width=1.5pt, mark size=4pt, mark=o, mark options={solid, mycolor4}]
  table[row sep=crcr]{%
64	7.903764\\
81	21.843459\\
100	59.875947\\
121	137.299153\\
};

\end{axis}

\end{tikzpicture}%
%
  \caption{Rate of convergence of the exponential Euler scheme
    \eqref{eq:expEullc}
    for the semidiscretization of
    ADR equation~\eqref{eq:ADR} with $n_1=n_2=n_3=n=20$ discretization points
    (left) and wall-clock time in seconds for increasing number $n$ of
    discretization points and 250 time steps
    up to final time $T=0.1$ (right).}
  \label{fig:ord1}
\end{figure}
After the semidiscretization in space by second-order centered
finite differences we end up with
an ODEs system of type~\eqref{eq:ODE}, with $K$ a matrix
with heptadiagonal structure.
This is a three-dimensional variation of a standard stiff example~\cite{HO05}
for exponential integrators.
\subsubsection{Exponential Euler}
We start by implementing the exponential Euler scheme \eqref{eq:expEullc} and
 test its correct order of convergence for a discretization
in space with $n_1=n_2=n_3=n=20$ internal points and a final simulation time
$T=0.1$ (Figure~\ref{fig:ord1}, left).
Then, we test the efficiency of the underlying methods for computing
linear combination~\eqref{eq:expEullc}, with a number of
discretization points in each direction
ranging from $n=64$ to $n=121$.
The results are presented in Figure~\ref{fig:ord1}, right.
We observe that \textsc{phiks} turns out to be at least twice as fast than the
other considered methods.
\begin{figure}[htb!]
  \centering
%
%
\definecolor{mycolor1}{rgb}{0.9922,0.9059,0.1451}%
\definecolor{mycolor2}{rgb}{0.3686,0.7882,0.38430}%
\definecolor{mycolor3}{rgb}{0.12940,0.56860,0.54900}%
\definecolor{mycolor4}{rgb}{0.23140,0.32160,0.54510}%
\definecolor{mycolor5}{rgb}{0.26670,0.00390,0.32940}%
\begin{tikzpicture}

\begin{axis}[%
width=1.8in,
height=1.5in,
at={(0.766in,0.486in)},
scale only axis,
xmode=log,
xmin=165,
xmax=480,
xminorticks=true,
xlabel style={font=\color{white!15!black}},
xlabel={number of time steps},
xtick={200,250,300,350,400},
xticklabels={{200},{250},{300},{350},{400}},
ymode=log,
ymin=1e-08,
ymax=1e-07,
yminorticks=true,
ylabel style={font=\color{white!15!black}},
ylabel={error},
axis background/.style={fill=white},
title style={font=\bfseries},
legend style={legend cell align=left, align=left, draw=white!15!black, font=\tiny},
legend image post style={scale=0.4}
]

\addplot [color=mycolor2, line width=1.5pt, only marks, mark size=4pt, mark=asterisk, mark options={solid, mycolor2}]
  table[row sep=crcr]{%
200	5.32365148375022e-08\\
250	3.40654373414008e-08\\
300	2.36537259976104e-08\\
350	1.73767359584327e-08\\
400	1.33032316313845e-08\\
};
\addlegendentry{\textsc{kiops}}

\addplot [color=mycolor3, line width=1.5pt, only marks, mark size=4pt, mark=+, mark options={solid, mycolor3}]
  table[row sep=crcr]{%
200	5.32364855058506e-08\\
250	3.40654173149628e-08\\
300	2.36537142649498e-08\\
350	1.73767355538582e-08\\
400	1.33032184827131e-08\\
};
\addlegendentry{\textsc{bamphi}}

\addplot [color=mycolor4, line width=1.5pt, only marks, mark size=4pt, mark=o, mark options={solid, mycolor4}]
  table[row sep=crcr]{%
200	5.32365767374016e-08\\
250	3.40654565586898e-08\\
300	2.36537361119731e-08\\
350	1.73767420270503e-08\\
400	1.33032346656933e-08\\
};
\addlegendentry{\textsc{phipm\_simul\_iom}}

\addplot [color=mycolor5, line width=1.5pt, only marks, mark size=4pt, mark=square, mark options={solid, mycolor5}]
  table[row sep=crcr]{%
200	5.32519766636632e-08\\
250	3.40717871382634e-08\\
300	2.36537310547918e-08\\
350	1.73767533551365e-08\\
400	1.33186400537108e-08\\
};
\addlegendentry{\textsc{phiks}}


\addplot [color=black, dashed, line width=1.5pt]
  table[row sep=crcr]{%
200	5.32129386627733e-08\\
400	1.33032346656933e-08\\
};
\addlegendentry{order two}

\end{axis}

\end{tikzpicture}%
\hfill
%
%
\definecolor{mycolor1}{rgb}{0.9922,0.9059,0.1451}%
\definecolor{mycolor2}{rgb}{0.3686,0.7882,0.38430}%
\definecolor{mycolor3}{rgb}{0.12940,0.56860,0.54900}%
\definecolor{mycolor4}{rgb}{0.23140,0.32160,0.54510}%
\definecolor{mycolor5}{rgb}{0.26670,0.00390,0.32940}%
\begin{tikzpicture}

\begin{axis}[%
width=1.8in,
height=1.5in,
at={(0.777in,0.482in)},
scale only axis,
xmin=54,
xmax=131,
xminorticks=true,
xlabel style={font=\color{white!15!black}},
xlabel={number $n$ of discretization points},
xtick={64,81,100,121},
xticklabels={{64},{81},{100},{121}},
ymode=log,
ymin=1,
ymax=1000,
yminorticks=true,
ylabel style={font=\color{white!15!black}},
ylabel={wall-clock time},
axis background/.style={fill=white},
title style={font=\bfseries},
legend style={legend cell align=left, align=left, draw=white!15!black, font=\tiny, at={(0.54,0.98)}},
legend image post style={scale=0.4}
]
\addplot [color=mycolor5, line width=1.5pt, mark size=4pt, mark=square, mark options={solid, mycolor5}]
  table[row sep=crcr]{%
64	3.97\\
81	9.91\\
100	2.15e1\\
121	4.95e1\\
};

\addplot [color=mycolor3, line width=1.5pt, mark size=4pt, mark=+, mark options={solid, mycolor3}]
  table[row sep=crcr]{%
64	1.35e1\\
81	2.70e1\\
100	9.63e1\\
121	2.15e2\\
};

\addplot [color=mycolor2, line width=1.5pt, mark size=4pt, mark=asterisk, mark options={solid, mycolor2}]
  table[row sep=crcr]{%
64	1.33e1\\
81	3.80e1\\
100	1.03e2\\
121	2.33e2\\
};

\addplot [color=mycolor4, line width=1.5pt, mark size=4pt, mark=o, mark options={solid, mycolor4}]
  table[row sep=crcr]{%
64	1.07e1\\
81	3.06e1\\
100	8.66e1\\
121	1.96e2\\
};

\end{axis}

\end{tikzpicture}%
%
  \caption{Rate of convergence of the ETD2RK scheme \eqref{eq:ETD2RKlc} for the
    semidiscretization of
    ADR
    equation~\eqref{eq:ADR} with $n_1=n_2=n_3=n=20$ discretization points
    (left) and wall-clock time in seconds for increasing number $n$ of
    discretization points and 100 time steps up to final
    time $T=0.1$ (right).}
  \label{fig:ord2}
\end{figure}
\subsubsection{Exponential Runge--Kutta scheme of order two}
We then move to the implementation of the ETD2RK scheme.
This integrator is implemented
following the linear combination approach {\color{black}(see
formula~\eqref{eq:ETD2RKlc}) and thus requires}
two calls of the algorithms
(see the beginning of Section~\ref{sec:lincomb}) for each time step.
Again, since we consider four
  different routines for the approximation of the actions of matrix functions,
we test the correct order of convergence of
the schemes. In addition, we measure
the performance of the routines
as the discretization is space
becomes finer and finer. The results, collected in
Figure~\ref{fig:ord2}, are similar to the exponential Euler case. We notice that
\textsc{phiks} turns
out to be almost four times faster than the best of the
other methods,
\textsc{phipm\_simul\_iom}, in the largest size scenario (total
number of degrees of freedom $N=121^3$).

\begin{table}[!ht]
\centering
  \caption{Wall-clock times in seconds and
    average number of Tucker operators {\color{black}$T_\#$}
    per time step 
    for the solution of the semidiscretized PDE~\eqref{eq:ADR} by
    ETD2RK implemented by two calls of \textsc{phiks} either for the linear
    combinations {\color{black}(see formula~\eqref{eq:ETD2RKlc}) or
    separately for the functions $\varphi_1$
    and $\varphi_2$ (see formula~\eqref{eq:ETD2RKdp}).}}
  \label{tab:ord2}
\begin{tabular}{c|cccc|cccc}
    &\multicolumn{4}{c|}{linear combination of $\varphi$-functions} &
    \multicolumn{4}{c}{$\varphi$-functions on the same vector} \\
    \hline
    $n$  & 64 & 81 & 100 & 121 & 64 & 81 & 100 & 121\\
wall-clock & 3.97 & 9.91 & 21.5 & 49.5 & 2.96 & 6.33 & 17.8 & 34.9 \\
{\color{black}$T_\#$} & {\color{black}20.0} & {\color{black}23.0} & {\color{black}23.0} 
& {\color{black}26.0} & {\color{black}12.0} & {\color{black}12.0} & {\color{black}15.0} & {\color{black}15.0}\\
\hline
  \end{tabular}
\end{table}

For comparison, we also implemented the integrator ETD2RK by two calls of
the routine \textsc{phiks} to
compute separately the actions of the functions $\varphi_1$ and $\varphi_2$
(see the beginning of Section~\ref{sec:diffphi}). The results
are reported in Table~\ref{tab:ord2} and show that this
approach leads to a smaller number
of Tucker operators which translates into less wall-clock time.
{\color{black}
  Predicting 
  which approach gives the smallest computational
  cost for a generic integrator is difficult.
  As a rule of thumb, we suggest to use the version of the algorithm
  which requires less calls and, when their number is the same, to prefer
  the computation of actions of $\varphi$-functions on the same vector,
  since the total number of 
  Tucker operators $T_\#$ is smaller (compare formulas~\eqref{eq:Tpsv}
  and~\eqref{eq:Tlcp}).}


\subsection{Allen--Cahn equation} \label{sec:ac}
In this section we
examine an example similar to the one reported in
Reference~\cite{MMPC22}, which describes the Sylvester approach for
the computation of the $\varphi$-functions. It is the
two-dimensional Allen--Cahn
phase-field model equation~\cite{FP03} for the concentration $u$
\begin{subequations}\label{eq:AC}
\begin{equation}
  \left\{
  \begin{aligned}
    \partial_t u(t,x_1,x_2) &= \Delta u(t,x_1,x_2)
    + \frac{1}{\epsilon^2}u(t,x_1,x_2)(1-u^2(t,x_1,x_2))\\
    &=\left(\Delta+\frac{1}{\epsilon^2}\right)u(t,x_1,x_2)+g(u(t,x_1,x_2)),\\
     u(0,x_1,x_2) &= u_0(x_1,x_2),
  \end{aligned}
  \right.
\end{equation}
  in the spatial domain $[0,1]^2$,
coupled with homogeneous Neumann boundary conditions.
The initial condition is given by
{\small
\begin{equation}
  u_0(x_1,x_2) = \tanh\left(
  \frac{\frac{1}{4} +\frac{1}{10}\cos\left(\beta\cdot\mathrm{atan2}
  \left(x_2-\frac{1}{2},x_1-\frac{1}{2}\right)\right)
  -\sqrt{\left(x_1-\frac{1}{2}\right)^2+\left(x_2-\frac{1}{2}\right)^2}}
  {\sqrt{2}\alpha}\right).
\end{equation}}%
\end{subequations}
We set $\epsilon=0.05$, $\beta=7$, $\alpha=0.75$,
and we discretize in space
with second-order centered finite differences, thus
obtaining a system in form~\eqref{eq:ODE} with $K$ a matrix with
pentadiagonal structure.
We simulate
until final time $T=0.025$.
Notice that the linear operator $\Delta+\frac{1}{\epsilon^2}$
guarantees a unique solution for the corresponding
Sylvester equation, even with homogeneous Neumann boundary
conditions.
As time marching scheme, we employ the third-order exponential Runge--Kutta
integrator with reduced tableau
\begin{equation}\label{eq:expRK3s3}
  \begin{array}{c|cc}
    c_2 &  & \\
    c_3 &\gamma c_2 \varphi_{2,2} + \frac{c_3^2}{c_2}\varphi_{2,3} & \\
    \hline
    & \frac{\gamma}{\gamma c_2+c_3}\varphi_2 & \frac{1}{\gamma c_2+c_3}\varphi_2
  \end{array}
\end{equation}
with $c_3=2c_2=1/2$ and
$\gamma=\frac{(3c_3-2)c_3}{(2-3c_2)c_2}=-4/5$,
see formula~(5.9) in Reference~\cite{HO05}.
{\color{black}Here and in the next tableaux $\varphi_{\ell,j}$ denotes
  $\varphi_\ell(c_j\tau K)$.} Its implementation
involves the usage only of the $\varphi_1$ and $\varphi_2$ functions,
which do not trigger the ill-conditioning of the Sylvester equation
observed in Reference~\cite{MMPC22}.
This integrator requires to compute
the following actions (scaled by proper coefficients)
\begin{equation*}
\begin{aligned}
  &\varphi_1(\tau K/2^j)\bb f(t_n,\bb u_n), & j&=0,1,2,\\
  &\varphi_2(\tau K/2^j)\bb d_{n2}, & j&=0,1,2,\\
  &\varphi_2(\tau K)\bb d_{n3},
\end{aligned}
\end{equation*}
see formula~\eqref{eq:RK}.
The \textsc{sylvphi} routine is then called in total six times: three
times to compute the action of the $\varphi_1$ function at the different
scales of $K$,
twice to compute the action of the $\varphi_2$ function for
the scales $j=1$ and $j=2$
and, finally, once to compute the action of $\varphi_2(\tau K)$ to
$(\gamma\bb d_{n2}+\bb d_{n3})/(\gamma c_2+c_3)$. Therefore,
nine Sylvester equations have to be solved. The other four routines
have to be called three times, one for each of the above
rows.
In fact, all of them are natively able to
produce the action of single $\varphi$-functions simultaneously
at different scales of $K$ (see,
in particular, Section~\ref{sec:diffphi} for \textsc{phiks}).
\begin{figure}[htb!]
  \centering
%
%
\definecolor{mycolor1}{rgb}{0.9922,0.9059,0.1451}%
\definecolor{mycolor2}{rgb}{0.3686,0.7882,0.38430}%
\definecolor{mycolor3}{rgb}{0.12940,0.56860,0.54900}%
\definecolor{mycolor4}{rgb}{0.23140,0.32160,0.54510}%
\definecolor{mycolor5}{rgb}{0.26670,0.00390,0.32940}%
\begin{tikzpicture}

\begin{axis}[%
width=1.8in,
height=1.5in,
at={(0.766in,0.486in)},
scale only axis,
xmode=log,
xmin=87,
xmax=230,
xminorticks=true,
xlabel style={font=\color{white!15!black}},
xlabel={number of time steps},
xtick={100,125,150,175,200},
xticklabels={{100},{125},{150},{175},{200}},
ymode=log,
ymin=1e-06,
ymax=1e-04,
yminorticks=true,
ylabel style={font=\color{white!15!black}},
ylabel={error},
axis background/.style={fill=white},
title style={font=\bfseries},
legend pos=south west,
legend style={legend cell align=left, align=left, draw=white!15!black, font=\tiny},
legend image post style={scale=0.4}
]

\addplot [color=mycolor1, line width=1.5pt, only marks, mark size=4.0pt, mark=triangle, mark options={solid, rotate=270, mycolor1}]
  table[row sep=crcr]{%
100	4.25883e-5\\
125	2.15302e-5\\
150	1.23586e-5\\
175	7.73981e-6\\
200	5.16689e-6\\
};
\addlegendentry{\textsc{sylvphi}}

\addplot [color=mycolor2, line width=1.5pt, only marks, mark size=4pt, mark=asterisk, mark options={solid, mycolor2}]
  table[row sep=crcr]{%
100	4.25882e-5\\
125	2.15301e-5\\
150	1.23587e-5\\
175	7.73982e-6\\
200	5.16690e-6\\
};
\addlegendentry{\textsc{kiops}}

\addplot [color=mycolor3, line width=1.5pt, only marks, mark size=4pt, mark=+, mark options={solid, mycolor3}]
  table[row sep=crcr]{%
100	4.25881e-5\\
125	2.15301e-5\\
150	1.23585e-5\\
175	7.73985e-6\\
200	5.16691e-6\\
};
\addlegendentry{\textsc{bamphi}}

\addplot [color=mycolor4, line width=1.5pt, only marks, mark size=4pt, mark=o, mark options={solid, mycolor4}]
  table[row sep=crcr]{%
100	4.25880e-5\\
125	2.15399e-5\\
150	1.23586e-5\\
175	7.73982e-6\\
200	5.16690e-6\\
};
\addlegendentry{\textsc{phipm\_simul\_iom}}

\addplot [color=mycolor5, line width=1.5pt, only marks, mark size=4pt, mark=square, mark options={solid, mycolor5}]
  table[row sep=crcr]{%
100	4.25881e-5\\
125	2.15398e-5\\
150	1.23585e-5\\
175	7.73980e-6\\
200	5.16691e-6\\
};
\addlegendentry{\textsc{phiks}}

\addplot [color=black, dashdotted, line width=1.5pt]
  table[row sep=crcr]{%
100	4.13351e-5\\
200	5.16689e-6\\
};
\addlegendentry{order three}

\end{axis}

\end{tikzpicture}%
\hfill
%
%
\definecolor{mycolor1}{rgb}{0.9922,0.9059,0.1451}%
\definecolor{mycolor2}{rgb}{0.3686,0.7882,0.38430}%
\definecolor{mycolor3}{rgb}{0.12940,0.56860,0.54900}%
\definecolor{mycolor4}{rgb}{0.23140,0.32160,0.54510}%
\definecolor{mycolor5}{rgb}{0.26670,0.00390,0.32940}%
\begin{tikzpicture}

\begin{axis}[%
width=1.8in,
height=1.5in,
at={(0.766in,0.486in)},
scale only axis,
xmin = 321,
xmax=681,
xminorticks=true,
xlabel style={font=\color{white!15!black}},
xlabel={number $n$ of discretization points},
xtick={351,451,551,651},
xticklabels={{351},{451},{551},{651}},
ymode=log,
ymin=1,
ymax=1e3,
yminorticks=true,
ylabel style={font=\color{white!15!black}},
ylabel={wall-clock time},
axis background/.style={fill=white},
title style={font=\bfseries},
legend style={legend cell align=left, align=left, draw=white!15!black, font=\footnotesize, at={(0.6,0.98)}},
legend image post style={scale=0.6}
]
\addplot [color=mycolor5, line width=1.5pt, mark size=4pt, mark=square, mark options={solid, mycolor5}]
  table[row sep=crcr]{%
351	7.71e0\\
451	1.82e1\\
551	3.48e1\\
651	5.17e1\\
};

\addplot [color=mycolor3, line width=1.5pt, mark size=4pt, mark=+, mark options={solid, mycolor3}]
  table[row sep=crcr]{%
351	1.68e1\\
451	4.34e1\\
551	8.01e1\\
651	1.42e2\\
};

\addplot [color=mycolor2, line width=1.5pt, mark size=4pt, mark=asterisk, mark options={solid, mycolor2}]
  table[row sep=crcr]{%
351	1.46e1\\
451	3.21e1\\
551	6.11e1\\
651	1.18e2\\
};

\addplot [color=mycolor4, line width=1.5pt, mark size=4pt, mark=o, mark options={solid, mycolor4}]
  table[row sep=crcr]{%
351	1.19e1\\
451	2.76e1\\
551	5.36e1\\
651	1.01e2\\
};

\addplot [color=mycolor1, line width=1.5pt, mark size=4pt, mark=triangle, mark options={solid, rotate=270, mycolor1}]
  table[row sep=crcr]{%
351	2.19e1\\
451	5.26e1\\
551	6.63e1\\
651	9.24e1\\
};

\end{axis}

\end{tikzpicture}%
%
  \caption{Rate of convergence of
    the exponential Runge--Kutta scheme~\eqref{eq:expRK3s3}
    for the semidiscretization of
    Allen--Cahn equation~\eqref{eq:AC} with $n_1=n_2=n=21$ discretization points
    (left) and wall-clock time in seconds for increasing number $n$ of
    discretization points and 20 time steps
    up to final time $T=0.025$ (right).}
  \label{fig:ord3}
\end{figure}
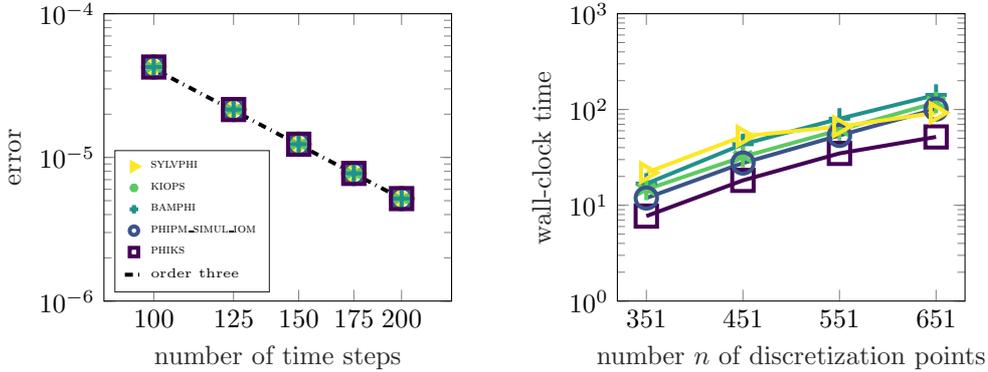
The results are summarized in Figure~\ref{fig:ord3}.
Also in this two-dimensional example, with numbers of degrees of freedom
up to $N=651^2$,
the \textsc{phiks} routine turns out to be always the fastest
by a factor of roughly 1.5 with respect to the other
techniques.
\subsection{Brusselator model}\label{sec:bruss}
In this section we apply the $\mu$-mode approach
to the block diagonal ODEs system which arises from the
semidiscretization in space of the three-dimensional Brusselator
model~\cite{GRT18,UW20} for the
two chemical concentrations $u$ and~$v$
\begin{equation}\label{eq:bruss}
  \left\{
  \begin{aligned}
    \partial_t u(t,x_1,x_2,x_3) &= d_1\Delta u(t,x_1,x_2,x_3)
    -(b+1)u(t,x_1,x_2,x_3)\\
    &+a+ u^2(t,x_1,x_2,x_3)v(t,x_1,x_2,x_3),\\
    \partial_t v(t,x_1,x_2,x_3) &= d_2\Delta v(t,x_1,x_2,x_3)\\
    &+ bu(t,x_1,x_2,x_3)
    -u^2(t,x_1,x_2,x_3)v(t,x_1,x_2,x_3),\\
    u(0,x_1,x_2,x_3) &= 64^2x_1^2(1-x_1)^2x_2^2(1-x_2)^2x_3^2(1-x_3)^2, \\
    v(0,x_1,x_2,x_3) &= c,
  \end{aligned}
  \right.
\end{equation}
in the spatial domain $[0,1]^3$. The system is
completed with homogeneous Neumann boundary conditions.
The semidiscretization in space by finite differences yields the system
\begin{equation}\label{eq:ODE2}
    \begin{pmatrix}
      \bb u'(t)\\
      \bb v'(t)\end{pmatrix}
    = \begin{pmatrix}
      K_1 & 0\\
      0 & K_2
    \end{pmatrix}
    \begin{pmatrix}
    \bb u(t)\\
    \bb v(t)
    \end{pmatrix}+
    \begin{pmatrix}
a+\bb u^2(t)\bb v(t)\\
b\bb u(t)- \bb u^2(t)\bb v(t)
    \end{pmatrix}\iff
    \bb w'(t)=\hat K \bb w(t)+\bb g(t,\bb w(t)),
\end{equation}
where $K_1$ is the discretization matrix of the linear
operator $d_1\Delta -(b+1)$ and $K_2$ is the discretization matrix
of the linear operator $d_2\Delta$, both clearly Kronecker
sums~\eqref{eq:kronsum}.
System~\eqref{eq:ODE2} cannot be written in Kronecker form~\eqref{eq:ODE}.
However, the actions of $\varphi$-functions can still be
efficiently computed by the proposed approach, since the
matrix $\hat K$ in system~\eqref{eq:ODE2} is block diagonal.
We set $a=c=1$, $b=3$, $d_1=d_2=0.02$ {\color{black}and, since
  we are going to employ two exponential integrators of order
  four, we also increase the accuracy in space by using} finite differences of
order four (leading to matrices $K_1$ and $K_2$ with a 13-diagonal structure),
and we simulate until final time $T=1$
employing two
fourth-order exponential integrators of Runge--Kutta type.
\subsubsection{An exponential Runge--Kutta scheme of order four with five
  stages}\label{sec:column-wise}
The first exponential integrator that we consider
has reduced tableau
\begin{equation}\label{eq:expRK4s5}
  \begin{array}{c|cccc}
    \frac{1}{2} &  & & &\\
    \frac{1}{2} & \varphi_{2,3}& & &\\
    1 & \varphi_{2,4} & \varphi_{2,4}& &\\
    \frac{1}{2} & a_{52} & a_{52} & \frac{1}{4}\varphi_{2,5}-a_{52}&\\
    \hline
    &   0 & 0 & -\varphi_2+4\varphi_3 & 4\varphi_2-8\varphi_3
  \end{array}
\end{equation}
with
\begin{equation*}
  a_{52}=\frac{1}{2}\varphi_{2,5}-\varphi_{3,4}+\frac{1}{4}\varphi_{2,4}
  -\frac{1}{2}\varphi_{3,5},
\end{equation*}
see tableau~(5.19) in Reference~\cite{HO05}.
It is possible to
effectively use the routine \textsc{phiks} to
evaluate
$\varphi$-functions applied to the same vector.
In this way, the following quantities have to be computed
\begin{equation*}
\begin{aligned}
  &\varphi_1(\tau \hat K/2^j)\bb f(t_n,\bb w_n),&&&j&=0,1,\\
  &  \varphi_\ell(\tau \hat K/2^j)\bb d_{n2},&\ell&=2,3,&j&=0,1,\\
  &\varphi_\ell(\tau \hat K/2^j)\bb d_{n3},&\ell&=2,3,&j&=0,1,\\
  &\varphi_\ell(\tau \hat K/2^j)\bb d_{n4},&\ell&=2,3,&j&=0,1,\\
  &\varphi_\ell(\tau \hat K)\bb d_{n5},&\ell&=2,3.
\end{aligned}
\end{equation*}
In each of the five lines, only two calls of
\textsc{phiks} are needed to compute the desired actions,
due to the block diagonal structure of the
matrix~$\hat K$. Indeed, in the first four lines
we can obtain
the actions of $\varphi$-functions
simultaneously at $j=0$ and $j=1$, thanks to the squaring
algorithm~\eqref{eq:simul}. In addition, in the second, third, and
fourth line,
\emph{both} $\varphi_2$ and $\varphi_3$ are produced at different scales.
This is not possible for the routines \textsc{phipm\_simul\_iom},
\textsc{kiops}, and \textsc{bamphi}. By proceeding with
a sequential implementation of the stages
in tableau~\eqref{eq:expRK4s5},
they would require six calls with the action of the
matrix $\hat K$, as already noticed in Reference~\cite{L21}. However,
we found that
it is possible to alternatively assemble the stages by computing the following
quantities
\begin{equation}\label{eq:expRK4s5_smart}
\begin{aligned}
  &\varphi_1(\tau \hat K/2^j)\bb f(t_n,\bb w_n),&j&=0,1,\\
  &\varphi_2(\tau \hat K/2^j)\bb d_{n2},&j&=0,1,\\
  &\varphi_2(\tau \hat K/2^j)\bb d_{n3},&j&=0,1,\\
  &\frac{\varphi_2(\tau \hat K/2^j)}{2^{2j}}\bb d_{n4}+
  \frac{\varphi_3(\tau \hat K/2^j)}{2^{3j}}(4\bb d_{n2}+4\bb d_{n3}-4\bb d_{n4}),&j&=0,1,\\
  &\varphi_2(\tau \hat K)(4\bb d_{n5}-\bb d_{n4})+
  \varphi_3(\tau \hat K)(4\bb d_{n4}-8\bb d_{n5}),
\end{aligned}
\end{equation}
which require again two calls to \textsc{phiks} for each
  of the five lines,
and in total five calls to the other routines.
To ensure a fair comparison, we follow
approach~\eqref{eq:expRK4s5_smart} with all the routines
under comparison.
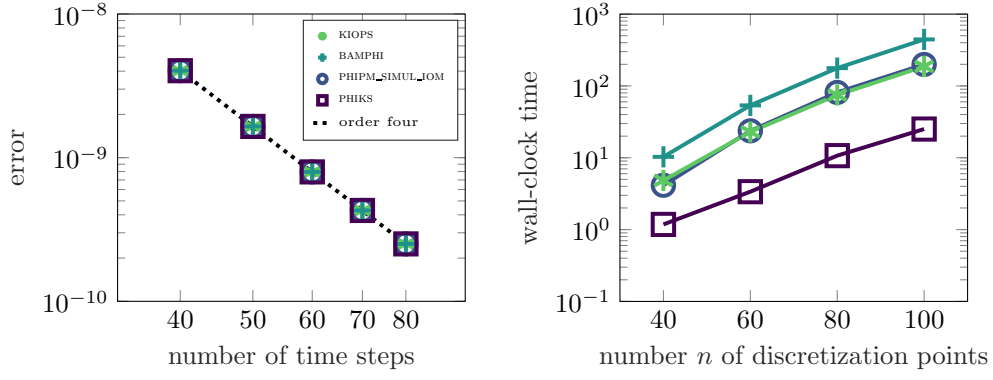
\begin{figure}[!htb]
  \centering
%
%
\definecolor{mycolor1}{rgb}{0.9922,0.9059,0.1451}%
\definecolor{mycolor2}{rgb}{0.3686,0.7882,0.38430}%
\definecolor{mycolor3}{rgb}{0.12940,0.56860,0.54900}%
\definecolor{mycolor4}{rgb}{0.23140,0.32160,0.54510}%
\definecolor{mycolor5}{rgb}{0.26670,0.00390,0.32940}%
\begin{tikzpicture}

\begin{axis}[%
width=1.8in,
height=1.5in,
at={(0.766in,0.486in)},
scale only axis,
xmode=log,
xmin=33,
xmax=96,
xminorticks=true,
xlabel style={font=\color{white!15!black}},
xlabel={number of time steps},
xtick={40,50,60,70,80},
xticklabels={{40},{50},{60},{70},{80}},
ymode=log,
ymin=1e-10,
ymax=1e-08,
yminorticks=true,
ylabel style={font=\color{white!15!black}},
ylabel={error},
axis background/.style={fill=white},
title style={font=\bfseries},
legend style={legend cell align=left, align=left, draw=white!15!black, font=\tiny},
legend image post style={scale=0.4}
]

\addplot [color=mycolor2, line width=1.5pt, only marks, mark size=4pt, mark=asterisk, mark options={solid, mycolor2}]
  table[row sep=crcr]{%
40	4.03e-09\\
50	1.65e-09\\
60	7.96e-10\\
70	4.29e-10\\
80	2.51e-10\\
};
\addlegendentry{\textsc{kiops}}

\addplot [color=mycolor3, line width=1.5pt, only marks, mark size=4pt, mark=+, mark options={solid, mycolor3}]
  table[row sep=crcr]{%
40	4.03e-09\\
50	1.65e-09\\
60	7.96e-10\\
70	4.29e-10\\
80	2.51e-10\\
};
\addlegendentry{\textsc{bamphi}}

\addplot [color=mycolor4, line width=1.5pt, only marks, mark size=4pt, mark=o, mark options={solid, mycolor4}]
  table[row sep=crcr]{%
40	4.03e-09\\
50	1.65e-09\\
60	7.96e-10\\
70	4.29e-10\\
80	2.51e-10\\
};
\addlegendentry{\textsc{phipm\_simul\_iom}}

  \addplot [color=mycolor5, line width=1.5pt, only marks, mark size=4.pt, mark=square, mark options={solid, mycolor5}]
  table[row sep=crcr]{%
40	4.03e-09\\
50	1.65e-09\\
60	7.96e-10\\
70	4.29e-10\\
80	2.51e-10\\
};
\addlegendentry{\textsc{phiks}}

\addplot [color=black, dotted, line width=1.5pt]
  table[row sep=crcr]{%
40	4.03e-09\\
80	2.51e-10\\
};
\addlegendentry{order four}

\end{axis}

\end{tikzpicture}%
\hfill
%
%
\definecolor{mycolor1}{rgb}{0.9922,0.9059,0.1451}%
\definecolor{mycolor2}{rgb}{0.3686,0.7882,0.38430}%
\definecolor{mycolor3}{rgb}{0.12940,0.56860,0.54900}%
\definecolor{mycolor4}{rgb}{0.23140,0.32160,0.54510}%
\definecolor{mycolor5}{rgb}{0.26670,0.00390,0.32940}%
\begin{tikzpicture}

\begin{axis}[%
width=1.8in,
height=1.5in,
at={(0.777in,0.482in)},
scale only axis,
xmin=30,
xmax=110,
xminorticks=true,
xlabel style={font=\color{white!15!black}},
xlabel={number $n$ of discretization points},
xtick={40,60,80,100},
xticklabels={{40},{60},{80},{100}},
ymode=log,
ymin=0.1,
ymax=1000,
yminorticks=true,
ylabel style={font=\color{white!15!black}},
ylabel={wall-clock time},
axis background/.style={fill=white},
title style={font=\bfseries},
legend style={legend cell align=left, align=left, draw=white!15!black, font=\footnotesize}
]

  \addplot [color=mycolor5, line width=1.5pt, mark size=4.0pt, mark=square, mark options={solid, mycolor5}]
  table[row sep=crcr]{%
40	1.18\\
60	3.37\\
80	1.07e1\\
100	2.51e1\\
};

  \addplot [color=mycolor4, line width=1.5pt, mark size=4.0pt, mark=o, mark options={solid, mycolor4}]
  table[row sep=crcr]{%
40	4.14\\
60	2.37e1\\
80	8.12e1\\
100	1.99e2\\
};

\addplot [color=mycolor3, line width=1.5pt, mark size=4.0pt, mark=+, mark options={solid, mycolor3}]
  table[row sep=crcr]{%
40	1.03e1\\
60	5.36e1\\
80	1.77e2\\
100	4.44e2\\
};

\addplot [color=mycolor2, line width=1.5pt, mark size=4.0pt, mark=asterisk, mark options={solid, mycolor2}]
  table[row sep=crcr]{%
40	4.83\\
60	2.29e1\\
80	7.48e1\\
100	1.87e2\\
};

\end{axis}

\end{tikzpicture}%
%
  \caption{Rate of convergence of the exponential Runge--Kutta
    scheme~\eqref{eq:expRK4s5}
    for the semidiscretization of
    Brusselator
    model~\eqref{eq:bruss} with $n_1=n_2=n_3=n=11$ discretization points
    (left) and wall-clock time in seconds for increasing number $n$ of
    discretization points and 20 time steps up to final time
    $T=1$ (right).}
\label{fig:ord4_1}
\end{figure}%
As we did in the previous numerical examples,
we check the correct order of convergence of the integrator for different
numbers of time steps and $n_1=n_2=n_3=n=20$,
and we measure the wall-clock time for different
numbers $n$
of discretization points for each dimension,
with $n$ ranging from 40 to 100, leading to a maximum number of degrees of
freedom $N=2\cdot 100^3$.
The integration is performed up to the final time $T=1$.
The results are collected in
Figure~\ref{fig:ord4_1}.
Again, \textsc{phiks} performs better than the other methods,
being up to 7.5 times faster.
The speed-up is larger than in the previous examples mainly due
to the choice of the spatial discretization, that
leads to denser matrices, and which affects all the routines but
\textsc{phiks} (as it is based on a $\mu$-mode approach).
The insensitivity of
the $\mu$-mode approach to the density of the matrices was
already pointed out in Reference~\cite{CCEOZ22}.
\subsubsection{An exponential Runge--Kutta scheme of order four with
  six stages}\label{sec:row-wise}
\begin{figure}[!ht]
  \centering
%
%
\definecolor{mycolor1}{rgb}{0.9922,0.9059,0.1451}%
\definecolor{mycolor2}{rgb}{0.3686,0.7882,0.38430}%
\definecolor{mycolor3}{rgb}{0.12940,0.56860,0.54900}%
\definecolor{mycolor4}{rgb}{0.23140,0.32160,0.54510}%
\definecolor{mycolor5}{rgb}{0.26670,0.00390,0.32940}%
\begin{tikzpicture}

\begin{axis}[%
width=1.8in,
height=1.5in,
at={(0.766in,0.486in)},
scale only axis,
xmode=log,
xmin=33,
xmax=96,
xminorticks=true,
xlabel style={font=\color{white!15!black}},
xlabel={number of time steps},
xtick={40,50,60,70,80},
xticklabels={{40},{50},{60},{70},{80}},
ymode=log,
ymin=1e-10,
ymax=1e-08,
yminorticks=true,
ylabel style={font=\color{white!15!black}},
ylabel={error},
axis background/.style={fill=white},
title style={font=\bfseries},
legend style={legend cell align=left, align=left, draw=white!15!black, font=\tiny},
legend image post style={scale=0.4}
]

\addplot [color=mycolor2, line width=1.5pt, only marks, mark size=4pt, mark=asterisk, mark options={solid, mycolor2}]
  table[row sep=crcr]{%
40	4.66e-09\\
50	1.915e-09\\
60	9.24e-10\\
70	4.98e-10\\
80	2.92e-10\\
};
\addlegendentry{\textsc{kiops}}

\addplot [color=mycolor3, line width=1.5pt, only marks, mark size=4pt, mark=+, mark options={solid, mycolor3}]
  table[row sep=crcr]{%
40	4.66e-09\\
50	1.915e-09\\
60	9.24e-10\\
70	4.98e-10\\
80	2.92e-10\\
};
\addlegendentry{\textsc{bamphi}}

\addplot [color=mycolor4, line width=1.5pt, only marks, mark size=4pt, mark=o, mark options={solid, mycolor4}]
  table[row sep=crcr]{%
40	4.68e-09\\
50	1.92e-09\\
60	9.25e-10\\
70	4.99e-10\\
80	2.92e-10\\
};
\addlegendentry{\textsc{phipm\_simul\_iom}}

\addplot [color=mycolor5, line width=1.5pt, only marks, mark size=4pt, mark=square, mark options={solid, mycolor5}]
  table[row sep=crcr]{%
40	4.67e-09\\
50	1.91e-09\\
60	9.24e-10\\
70	4.98e-10\\
80	2.92e-10\\
};
\addlegendentry{\textsc{phiks}}

\addplot [color=black, dotted, line width=1.5pt]
  table[row sep=crcr]{%
40	4.67e-09\\
80	2.92e-10\\
};
\addlegendentry{order four}

\end{axis}

\end{tikzpicture}%
\hfill
%
%
\definecolor{mycolor1}{rgb}{0.9922,0.9059,0.1451}%
\definecolor{mycolor2}{rgb}{0.3686,0.7882,0.38430}%
\definecolor{mycolor3}{rgb}{0.12940,0.56860,0.54900}%
\definecolor{mycolor4}{rgb}{0.23140,0.32160,0.54510}%
\definecolor{mycolor5}{rgb}{0.26670,0.00390,0.32940}%
\begin{tikzpicture}

\begin{axis}[%
width=1.8in,
height=1.5in,
at={(0.766in,0.486in)},
scale only axis,
xmin=30,
xmax=110,
xminorticks=true,
xlabel style={font=\color{white!15!black}},
xlabel={number $n$ of discretization points},
xtick={40,60,80,100},
xticklabels={{40},{60},{80},{100}},
ymode=log,
ymin=0.1,
ymax=1000,
yminorticks=true,
ylabel style={font=\color{white!15!black}},
ylabel={wall-clock time},
axis background/.style={fill=white},
title style={font=\bfseries},
legend style={legend cell align=left, align=left, draw=white!15!black}
]
\addplot [color=mycolor5, line width=1.5pt, mark size=4.0pt, mark=square, mark options={solid, mycolor5}]
  table[row sep=crcr]{%
40	8.17e-1\\
60	2.44e0\\
80	8.22e0\\
100	2.14e1\\
};

\addplot [color=mycolor4, line width=1.5pt, mark size=4.0pt, mark=o, mark options={solid, mycolor4}]
  table[row sep=crcr]{%
40	3.73e0\\
60	2.01e1\\
80	7.13e1\\
100	1.79e2\\
};

\addplot [color=mycolor3, line width=1.5pt, mark size=4.0pt, mark=+, mark options={solid, mycolor3}]
  table[row sep=crcr]{%
40	6.51e0\\
60	3.14e1\\
80	9.17e1\\
100	2.49e2\\
};

\addplot [color=mycolor2, line width=1.5pt, mark size=4.0pt, mark=asterisk, mark options={solid, mycolor2}]
  table[row sep=crcr]{%
40	3.90e0\\
60	2.15e1\\
80	7.08e1\\
100	1.80e2\\
};

\end{axis}

\end{tikzpicture}%
%
  \caption{Rate of convergence of the exponential Runge--Kutta
    scheme~\eqref{eq:expRK4s6}
    for the semidiscretization of
    Brusselator
    model~\eqref{eq:bruss} with $n_1=n_2=n_3=n=11$ discretization points
    (left) and wall-clock time in seconds for increasing number $n$ of
    discretization points and 20 time steps up to final time
    $T=1$ (right).}
\label{fig:ord4_2}
\end{figure}
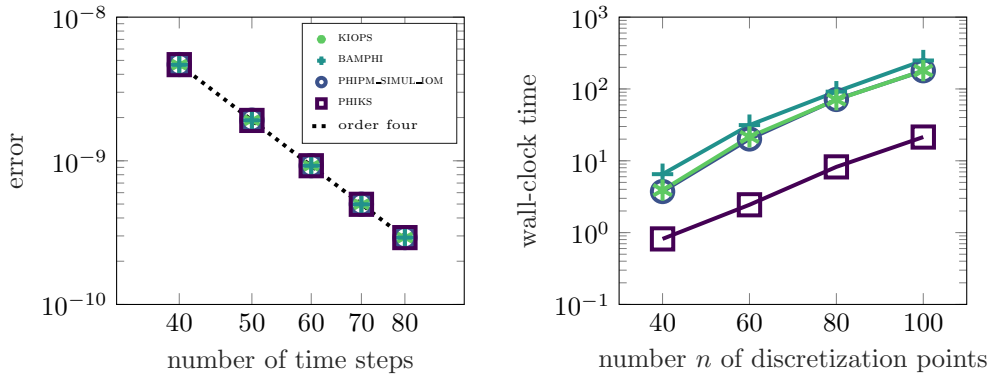

Finally, we consider the exponential Runge--Kutta integrator
of order four with reduced tableau
\begin{equation}\label{eq:expRK4s6}
  \begin{array}{c|ccccc}
    c_2 &  & & & &\\
    c_3 &  \frac{c_3^2}{c_2}\varphi_{2,3}& & & &\\
    c_4 &  \frac{c_4^2}{c_2}\varphi_{2,4} &  & & &\\
    c_5 &  0 & a_{5,3} & a_{5,4}
    &\\
    c_6 &   0& a_{6,3}&
    a_{6,4}\\
    \hline
    &  0 & 0 & 0 & b_5& b_6
  \end{array}
\end{equation}
where
\begin{equation*}
  \begin{aligned}
a_{5,3}&=\frac{c_4c_5^2}{c_3(c_4-c_3)}\varphi_{2,5}+
\frac{2c_5^3}{c_3(c_3-c_4)}\varphi_{3,5},&
a_{5,4}&=\frac{c_3c_5^2}{c_4(c_3-c_4)}\varphi_{2,5}+
\frac{2c_5^3}{c_4(c_4-c_3)}\varphi_{3,5},\\
a_{6,3}&=\frac{c_4c_6^2}{c_3(c_4-c_3)}\varphi_{2,6}+
\frac{2c_6^3}{c_3(c_3-c_4)}\varphi_{3,6},&
a_{6,4}&=\frac{c_3c_6^2}{c_4(c_3-c_4)}\varphi_{2,6}+
\frac{2c_6^3}{c_4(c_4-c_3)}\varphi_{3,6},\\
b_5&=\frac{c_{6}}{c_5(c_{6}-c_5)}\varphi_2
    +\frac{2}{c_5(c_5-c_{6})}\varphi_3,&
b_6&=\frac{c_{5}}{c_6(c_{5}-c_6)}\varphi_2
    +\frac{2}{c_6(c_6-c_{5})}\varphi_3,
  \end{aligned}
  \end{equation*}
and $c_3\ne c_4$, $c_5\ne c_6$, $c_6=(4c_5-3)/(6c_5-4)$,
see equation~(4.10) in Reference~\cite{L21}.
This integrator was designed so that multiple stages can be
computed simultaneously. In fact,
stages three and four and stages five and six require the same
combination of $\varphi$-functions
at different intermediate times, and therefore they can
be efficiently implemented by the routines \textsc{kiops},
\textsc{bamphi}, and \textsc{phipm\_simul\_iom} (the one
which has been originally employed with this integrator).
Here, by selecting $c_4=2c_3=2/3$
and $c_6=2c_5=1$ we can do the same with
the routine \textsc{phiks},
since the evaluation of a linear combination at a half time
comes for free by the new squaring algorithm~\eqref{eq:comb}.
The results, by setting the remaining free coefficient $c_2=1/3$,
are reported in Figure~\ref{fig:ord4_2}.
Since the number of calls needed by this integrator is smaller than
in the previous case, all the methods turn out to be slightly faster.
The speed-up of \textsc{phiks} with respect to the other routines
ranges from 4.5 to 8.6, depending on the size of the problem.
\section{Conclusions}\label{sec:conclusions}
In this manuscript, we proposed an efficient $\mu$-mode approach
to compute actions of $\varphi$-functions for matrices
$K$ which are
Kronecker sums of any number of arbitrary matrices $A_\mu$.
This structure naturally arises when suitably discretizing in space
some evolutionary PDEs of great importance in science and engineering,
such as advection--diffusion--reaction, Allen--Cahn, or Brusselator
equations, among the others. The corresponding stiff system of ODEs
can be effectively solved by
exponential integrators, which
rely on the efficient approximation of the action
of single $\varphi$-functions or linear combinations of them.
Our new method, that we named \textsc{phiks},
approximates the integral definition
of $\varphi$-functions by the Gauss--Lobatto--Legendre quadrature formula,
employs scaling and squaring techniques,
and computes the required actions in a
$\mu$-mode fashion by means of Tucker operators and exponentials
of the small sized matrices $A_\mu$,
exploiting the efficiency of modern hardware architectures to perform
level 3 BLAS operations.
Moreover, it allows to compute the desired quantities at
  different time scales, feature of great importance in the context
  of high order exponential integrators.
We tested our approach on different
stiff ODEs systems arising from the discretization of important
PDEs in two and three space dimensions,
using different exponential integrators (from stiff
order one to four) and different discretization matrices (finite
differences of order two or four).
As terms of comparison, we considered another
technique for computing actions of $\varphi$-functions
of Kronecker sums of matrices (based on Sylvester
equations, and currently limited to two space dimensions)
and more general techniques for computing actions of
$\varphi$-functions. Our method turned out to be always faster than
the others, with speed-ups ranging from 1.5 to 8.6, depending on the
example under consideration.
{\color{black}The numerical experiments
  suggest that \textsc{phiks} is preferable to existing methods,
  in particular for problems with denser matrices and
  for exponential integrators that can be implemented by computing
  suitable scales of the
  underlying (linear combinations of) $\varphi$-functions.}
Interesting future developments are 
the application of the method to space-fractional diffusion
equations~\cite{ZWC18}
and second-order in time partial differential equations~\cite{PO22}.
\section*{Declarations}
\subsection*{Acknowledgments}
Marco Caliari and Fabio Cassini are members of the Gruppo Nazionale
Calcolo Scientifico-Istituto Nazionale di Alta Matematica (GNCS-INdAM).
\subsection*{Funding}
The authors received partial support from the
Program Ricerca di Base 2019 of the University of Verona entitled
``Geometric Evolution of Multi Agent Systems''.
Fabio Cassini received financial support from
the Italian Ministry of University and Research (MUR) with the
PRIN Project 2022 No.~2022N9BM3N
``Efficient numerical schemes and optimal control methods for
time-dependent partial differential equations''.
\subsection*{Conflict of interest}
The authors declare no potential conflict of interests.
\bibliography{phiksbiblio}
\end{document}